\documentclass[10pt]{article}
\usepackage{graphicx}
\usepackage{amssymb}
\usepackage{amsmath}
\usepackage{epsfig}
\usepackage[english]{babel}
\usepackage{color} 

\def\div{{\rm div }}
\def\L{\mathcal{L}}
\def\tM{\tilde{M}}
\def\tR{R}

\def\R{\mathbb{R}}
\def\T{\mathbb{T}}
\def\Z{\mathbb{Z}}

\def\W{\mathcal{W}}
\def\D{\mathcal{D}}
\def\e{\mathcal{E}}
\def\V{\mathcal{V}}
\def\Y{Y}
\def\lameps{\lambda_{\varepsilon}}

\def\M{\mathcal{M}}
\def\psixi{\psi^{\xi,I}}
\def\psixiinf{\psixi_{\infty}}
\def\tpsiz{\tilde{\psi}_0}
\def\tpsio{\tilde{\psi}_1}
\def\tpsizi{\tilde{\psi}_i}
\def\tpsizinf{\tilde{\psi}_{\infty,0}}

\def\tpsioinf{\tilde{\psi}_{\infty,1}}
\def\tpsizinfi{\tilde{\psi}_{\infty}}

\def\sumi{\sum_{i=0}^1}

\def\psixibar{{\psi}^{\xi}}
\def\ffi{\frac{\psi^{\xi}}{\psixiinf}}

\def\ln{\text{ln}}
\def\max{\text{max}}
\def\exp{\text{exp}}
\def\E{\mathbb{E}} 

\newtheorem{theorem}{Theorem}
\newtheorem{definition}{Definition}
\newtheorem{lemma}{Lemma}
\newtheorem{remark}{Remark}

\def\proof{ {\it Proof : }}
\def\endproof{\hfill$\diamondsuit$\par\medskip}

\title{Long-time convergence of an Adaptive Biasing Force method: the bi-channel case}
\author{T. Leli\`evre and K. Minoukadeh \\ \\
\footnotesize {\'Ecole des Ponts ParisTech, CERMICS,\'Ecole des Ponts ParisTech, 6-8 avenue Blaise-Pascal,} \\
\footnotesize{77455 Champs-sur-Marne, Marne-la-Vallée cedex 2, France} \\
\footnotesize{MICMAC Project-Team, INRIA Rocquencourt, 78153 Le Chesnay, France}
}

\date{}
\begin{document}

\maketitle
\begin{abstract}
We present convergence results for an adaptive algorithm to compute free energies, namely the adaptive biasing force (ABF) method~\cite{darve-pohorille-01,henin2004}.  The free energy is the effective potential associated to a so-called reaction coordinate  $\xi(q)$, where $q = (q_1,\ldots,q_{3N})$ is the position vector of a $N$-particle system.  Computing free energy differences remains an important challenge in molecular dynamics due to the presence of meta-stable regions in the potential energy surface.  The ABF method uses an on-the-fly estimate of the free energy to bias dynamics and overcome metastability.  Using entropy arguments and logarithmic Sobolev inequalities, previous results have shown that the rate of convergence of the ABF method is limited by the metastable features of the canonical measures conditioned to being at fixed values of~$\xi$~\cite{lrs_nonlinearity}.  In this paper, we present an improvement on the existing results, in the presence of such metastabilities, which is a generic case encountered in practice.  More precisely, we study the so-called bi-channel case, where two channels along the reaction coordinate direction exist between an initial and final state, the channels being separated from each other by a region of very low probability.  With hypotheses made on `channel-dependent' conditional measures, we show on a bi-channel model that we introduce, that the convergence of the ABF method is in fact not limited by metastabilities in directions orthogonal to $\xi$ under two crucial assumptions: (i) exchange between the two channels is possible for {\em some} values of $\xi$ and (ii) the free energy is a good bias in each channel.  This theoretical result supports recent numerical experiments~\cite{minoukadeh2010}, where the efficiency of the ABF approach is demonstrated for such a multiple-channel situation.  
\end{abstract}

\section{Introduction}
\label{sec:introduction}
We consider a system of $N$ particles with positions $q \in \D \subset \R^{3N}$. In statistical physics, one is interested in calculating averages with respect to the Boltzmann-Gibbs measure
\begin{equation}
\label{eq:langevin}
 d\nu(q) = Z^{-1}\exp(-\beta V(q))dq, 
\end{equation}
\noindent with $V:\D \rightarrow \R$ the potential energy function, $Z = \int_{\D}\exp(-\beta V(q)) \ dq$ the partition function, $\D = \{q, V(q) < \infty\}$ the configuration space, and $\beta = 1/ (k_BT)$, where $k_B$ is the Boltzmann constant and $T$ is the temperature.  The function $V$ is the energy associated with the positions of the particles and is assumed to be such that $Z< \infty$.   The probability measure $\nu$ represents the equilibrium measure sampled by the particles in the canonical ensemble.  A typical dynamics that can be used to sample this measure through trajectorial averages is the overdamped Langevin dynamics
\begin{equation}\label{eq:sde} dQ_t = -\nabla V(Q_t)dt + \sqrt{2\beta^{-1}}dB_t,\end{equation}
where $(B_t)_{t \ge 0}$ is a $3N$-dimensional standard Brownian motion. Indeed, under loose assumptions on $V$, one has the ergodic property: for any smooth test function $\varphi$,
$$\lim_{T \to \infty} \frac{1}{T}\int_0^T \varphi(Q_t) \, dt = \int_{\D} \varphi \, d \nu.$$

The efficiency of this sampling procedure, which can be shown to be related to the convergence rate to equilibrium of the above dynamics is often hindered by {\em metastabilities in the potential function $V$}, namely regions of low energy are separated by high energy barriers. Equivalently, in terms of the probability measure, $\nu$ is typically a {\em multimodal measure}, with regions of high probability separated by regions of low probability. To circumvent this issue, a one-dimensional collective variable (or reaction coordinate) $\xi:\D\rightarrow\M$ is introduced, which will be used to define a biasing potential for~\eqref{eq:sde}.  In the following, we will assume that $|\nabla \xi| >0$ on $\D$, and that $\M = \T$, where $\T = \R/\Z$ is the one-dimensional torus (which typically corresponds to the case where the reaction coordinate represents a dihedral angle, for example to characterize the conformation of a molecule). Before defining more precisely the biased dynamics in the next section, we need to introduce a few notation.

 The image of the measure $\nu$ in $\xi$ is given by 
\begin{equation*}
 d\nu^{\xi}(z) = Z^{-1}\text{exp}(-\beta A(z))\ dz,
\end{equation*}
\noindent where $A$ is the so-called \textit{free energy}, defined by
 \begin{equation}
  \label{eq:A}
 A(z) = -\beta^{-1}\ln(Z_z)
 \end{equation}
where
\begin{equation}\label{eq:Z_z}
Z_z = \int_{\Sigma_z}\exp(-\beta V(q))\ \delta_{\xi(q)-z}(dq)
\end{equation}
is the partition function on the submanifold $\Sigma_z = \{q \in \D \ |\ \xi(q) = z\}$.  The measure $\delta_{\xi(q)-z}(dq)$ is defined through the conditioning formula: for any smooth test function $\varphi:\D \rightarrow \R$,
$$\int_{\D}\varphi(q)~dq = \int_{\M}\int_{\Sigma_z}\varphi(q) \delta_{\xi(q)-z}(dq)~dz.$$
Using the co-area formula (see~\cite[Appendix A]{lrs_nonlinearity}), one can also identify this measure as
 $\delta_{\xi(q)-z}(dq)=|\nabla \xi|^{-1}d\sigma_{\Sigma_z}$, where $\sigma_{\Sigma_z}$ is the Lebesgue measure on $\Sigma_z$. We assume in the following that $\xi$ and $V$ are such that $Z_z < \infty$, for all $z \in \M$.  

Practitioners are typically interested in free energy differences $A(z)-A(z_0)$, which can be obtained by computing (and integrating) the derivative $A'(z)$, the so-called {\em mean force} 
\begin{equation}
 \label{eq:ap}
A'(z) = Z_z^{-1}\int_{\Sigma_z}f(q)\ \exp(-\beta V(q))\ \delta_{\xi(q)-z}(dq)
\end{equation}
where $f$ is the \textit{local mean force} defined by
\begin{equation}
 \label{eq:F}
f =  \frac{\nabla V \cdot \nabla \xi}{|\nabla \xi|^2} 
           - \beta^{-1} \div \left(\frac{\nabla \xi}{|\nabla \xi|^2}\right).
\end{equation}
\noindent   The function $f$ can be understood as the negative force projected onto the reaction coordinate, plus some correction term related to the curvature of the submanifolds $\Sigma_z$.  Notice that the mean force~\eqref{eq:ap} is in fact a conditional expectation
\begin{equation}
 \label{eq:ap_exp} \displaystyle
A'(z) = \E_{\nu}[f(Q)|\xi(Q)=z] = \int_{\Sigma_z} f d \nu_{|z},  
\end{equation}
where $$d\nu_{|z}=\frac{\exp(-\beta V(q)) \delta_{\xi(q)-z}(dq)}{\exp(-\beta(A(z)))}$$
denotes the probability measure $\nu$ conditioned to a fixed value $z$ of $\xi(q)$. This measure is supported on the submanifold $\Sigma_z$. For the derivation of \eqref{eq:ap}--\eqref{eq:F}--\eqref{eq:ap_exp} which is again based on the co-area formula, the reader is referred to~\cite{ciccotti2008,sprik1998,otter1998}.

\subsection{The adaptive biasing force method}
The adaptive biasing force (ABF) method~\cite{darve-pohorille-01,henin2004} uses an estimate of the mean force $A'$ to bias the standard overdamped Langevin dynamics~\eqref{eq:sde} in order to overcome metastabilities in $\xi$. The bottom line of the approach is that it should be easier to sample the probability measure with density proportional to $\exp (- \beta ( V - A \circ \xi ) )$ than the original Boltzmann-Gibbs measure $\nu$, since the marginal probability of the former in $\xi$ is a uniform probability measure on $\T$, while the marginal of the latter (namely $\exp(-\beta A(z)) \, dz$) is typically multimodal.

The ABF dynamics is given by
\begin{equation}
\label{eq:abf_orig}
\left\{
\begin{aligned}
 dX_t     &= -\nabla(V-A_t \circ \xi)(X_t)dt + \sqrt{2\beta^{-1}}dB_t, \\
 A'_t (z) &= \E\left[f(X_t)|\xi(X_t)=z\right],
\end{aligned}
\right.
\end{equation}
\noindent where $A'_t$ is an on-the-fly estimate of the mean force, which, in view of the definition \eqref{eq:ap_exp}, is expected to be a good estimate of $A'$.  The law of $X_t$ has density $\psi(t,\cdot)$, which satisfies the non-linear Fokker-Planck equation:
\begin{equation}
\label{eq:pde_orig}
\left\{
\begin{array}{l}
 \partial_t \psi = \div(\nabla(V-A_t\circ \xi)\psi) + \beta^{-1} \Delta \psi,\\[5pt]
A_t'(z) = \frac{\displaystyle  \int_{\Sigma_z} f \ \psi \ \delta_{\xi(q)-z}(dq)}
		 {\displaystyle \int_{\Sigma_z}\psi\  \delta_{\xi(q)-z}(dq)}.
\end{array}
\right.
\end{equation}
\noindent Roughly speaking, the biasing force $\nabla (A_t \circ \xi)$ ``flattens the free-energy barriers in $\xi$''. To support this claim, let us simply indicate that if $|\nabla \xi|$ is constant,  the marginal density in $\xi$ satisfies a simple heat equation, with zero bias, see~\cite{lrs2007,lrs_nonlinearity} and also Equation~\eqref{eq:FPpsixibar} below.  Existence and uniqueness of solutions to  \eqref{eq:abf_orig} are studied in~\cite{jlr_2009} and a study of the longtime convergence of~\eqref{eq:pde_orig} can be found in~\cite{lrs_nonlinearity}, the results of which are briefly discussed below.

\subsection{Existing convergence results, and the multiple channel scenario}
It has been shown in~\cite{lrs_nonlinearity} that, under appropriate assumptions, the biasing force $A_t'$ in~\eqref{eq:abf_orig} (actually for a slightly different dynamics which reduces to~\eqref{eq:abf_orig} if $|\nabla \xi|$ is constant for example) converges to the mean force $A'$ exponentially fast in the longtime limit.  The rate of convergence was estimated as the minimum of (i) the rate at which the law of $\xi(X_t)$ converges to equilibrium, and (ii) the smallest logarithmic Sobolev inequality constant (LSI constant, discussed in Section~\ref{sec:entropy_defs}) of the conditional probability measures $\nu_{|z}$, for $z \in \M$. Thanks to the bias in the direction of the reaction coordinate, it can be shown that the rate of convergence of the marginal in $\xi$ is actually not the limiting rate in practice since it satisfies a simple diffusion equation. The real limitation is thus the metastable features ({\em i.e.} the multimodality) of the family of laws $\nu_{|z}$, which is quantified through the LSI constants associated to these measures: roughly speaking, the smaller the constant the more multimodal the probability measure. These constants may be in some cases extremely small, at least for some values of $z \in \M$. The question we address in this paper is the {\em optimality} of this theoretical rate of convergence for the ABF method.

The generic situation is indeed that the LSI constants for the measures $\nu_{|z}$ are not small uniformly in $z$. Typically, there exists some values of $z$ for which these measures are multimodal. This situation is  often encountered in practice due to two reasons.  First, finding a suitable reaction coordinate, namely in our context one that ensures that there are no metastabilities associated to the equilibrium measures $\nu_{|z}$, is not trivial for a large-dimensional system.  Secondly, a low-dimensional reaction coordinate may simply not be sufficient to describe all metastabilities of the system. In such cases, the results of~\cite{lrs_nonlinearity} predict a very small rate of convergence for the ABF dynamics~\eqref{eq:abf_orig}, and thus, the inefficiency of this biasing procedure.

As a typical case for which such difficulties appear, we will consider in the following the so-called `multiple channel situation' (see  Figure~\ref{fig:bichannel} on the left, for a bi-channel case). In such a situation, starting from a metastable state and as the system evolves in the direction of increasing values of the reaction coordinate $\xi$, it can follow different channels, which are separated (in the `orthogonal direction to $\xi$') by arbitrary high energy barriers.  In other words, the energy landscape features parallel valleys which are orthogonal to the isocontours of the reaction coordinate. In such a prototypical situation, the conditional probability measures $\nu_{|z}$ are indeed multimodal, for the values of $z$ corresponding to the system being in one of these channels.

However, recent numerical experiments~\cite{minoukadeh2010} (based on a discretization of the ABF dynamics~\eqref{eq:abf_orig} by multiple walkers simulated in parallel) suggest that in fact, high energy barriers in $\Sigma_z$ do not always hinder the convergence of the ABF method. The multiple walkers are made to follow similar dynamics~\eqref{eq:abf_orig}, but driven by independent Brownian motions. The chemical system considered in~\cite{minoukadeh2010} is the compact states of the deca-alanine peptide (the reaction coordinate $\xi$ being thus the end-to-end distance of the peptide). Due to some `buckling effects', this is a typical multiple channel situation, since the molecule can shrink to various compact states  (see~\cite{henin2010}). In~\cite{minoukadeh2010}, numerical results show that the ABF approach indeed yields reliable results in such a situation. We interpret this as follows. When encountered with a fork in the channel, each walker is likely to travel down a different channel. Thus, it is indeed almost impossible for a given walker to switch from one channel to another, once it has entered one of them, but this appears not to be necessary to obtain reliable results. It suggests that the theoretical rate of convergence obtained in~\cite{lrs_nonlinearity} is actually not optimal.

Inspired by these numerical results, we present herein an improved theoretical rate of convergence of the ABF method.  The rate will be shown to depend on the LSI constants of the family of equilibrium measures conditioned to being in $\Sigma_z$ {\em and} a channel. By doing so, we show that high energy barriers separating the channels do not in fact affect the rate of convergence of the method. The crucial assumptions needed to show our result are: (i) exchange between the two channels is possible for {\em some} values of $\xi$ (see [H1] below) and (ii) the free energy is a `good bias' in each channel (see [H4] below). This is formalized in the main result of this paper, namely Theorem~\ref{th:main} below.

For some technical reasons, we were actually unable to prove this result on the original ABF dynamics~\eqref{eq:abf_orig}. We will thus consider a slightly different system (that we call {\em the bi-channel model}) which retains the most important features of the dynamics~\eqref{eq:abf_orig} when applied to a potential exhibiting two parallel channels in the direction of $\xi$, separated by a high energy barrier. 

The paper is organized as follows.  In Section \ref{sec:model} we give details of the bi-channel model, define some probability measures and recall some entropy definitions before presenting the main result.  Finally, the proof of the main result is given in Section \ref{sec:proof}.

\section{The bi-channel model and statement of the main results}
\label{sec:model}
In this section, we present a model to describe the bi-channel scenario.  In the following, we treat the case $d=2$ (so that the position is $q=(x,y)$), $\D = \T \times \R$, $\M = \T$ and $\xi: \D \rightarrow \T$, where $\xi(x,y) = x$.  We further assume without loss of generality that $\beta = 1$, which can be done by a change of variables $\tilde{t} = \beta^{-1}t$, $\tilde{\psi}(\tilde{t},x) = \psi(t,x)$, $\tilde{V}(x,y) = \beta V(x,y)$.  With these assumptions, some notation may be simplified: $|\nabla \xi| = 1$, $\Sigma_z = \R$ and $\delta_{\xi(q)-z}(dq) = dy$.  Furthermore, the definition of the local mean force in \eqref{eq:F} simplifies to 
$$f=\partial_xV $$
and the free energy and its derivative are given by
\begin{equation}
\label{eq:act_A'}
A(x) = -\ln \int_{\R} \exp(- V(x,y))~dy 
 ~~~~ \text{and} ~~~~ 
A'(x) = \frac{\displaystyle \int_{\R} \partial_xV(x,y)~\exp(-V(x,y))~dy}{\displaystyle \int_{\R} \exp(- V(x,y))~dy}. 
\end{equation}
We would like to emphasize that the choice of the domain $\M = \T$ and reaction coordinate $\xi(x,y)~=~x$ is merely to reduce technicalities, see~\cite{lrs_nonlinearity} for appropriate tools to treat general~$\xi$.  In particular, the results can be straightforwardly generalized to the case $\D = \T \times \R^{d-1}$ and $\xi(q_1, \ldots,q_d)=q_1$. Likewise, the generalization to a situation with multiple channels (more than two) is straightforward (see Remark~\ref{rem:multiple_channels} below for another generalization).

\subsection{The bi-channel model}
\begin{figure}
\centering
\includegraphics[width=0.8\textwidth]{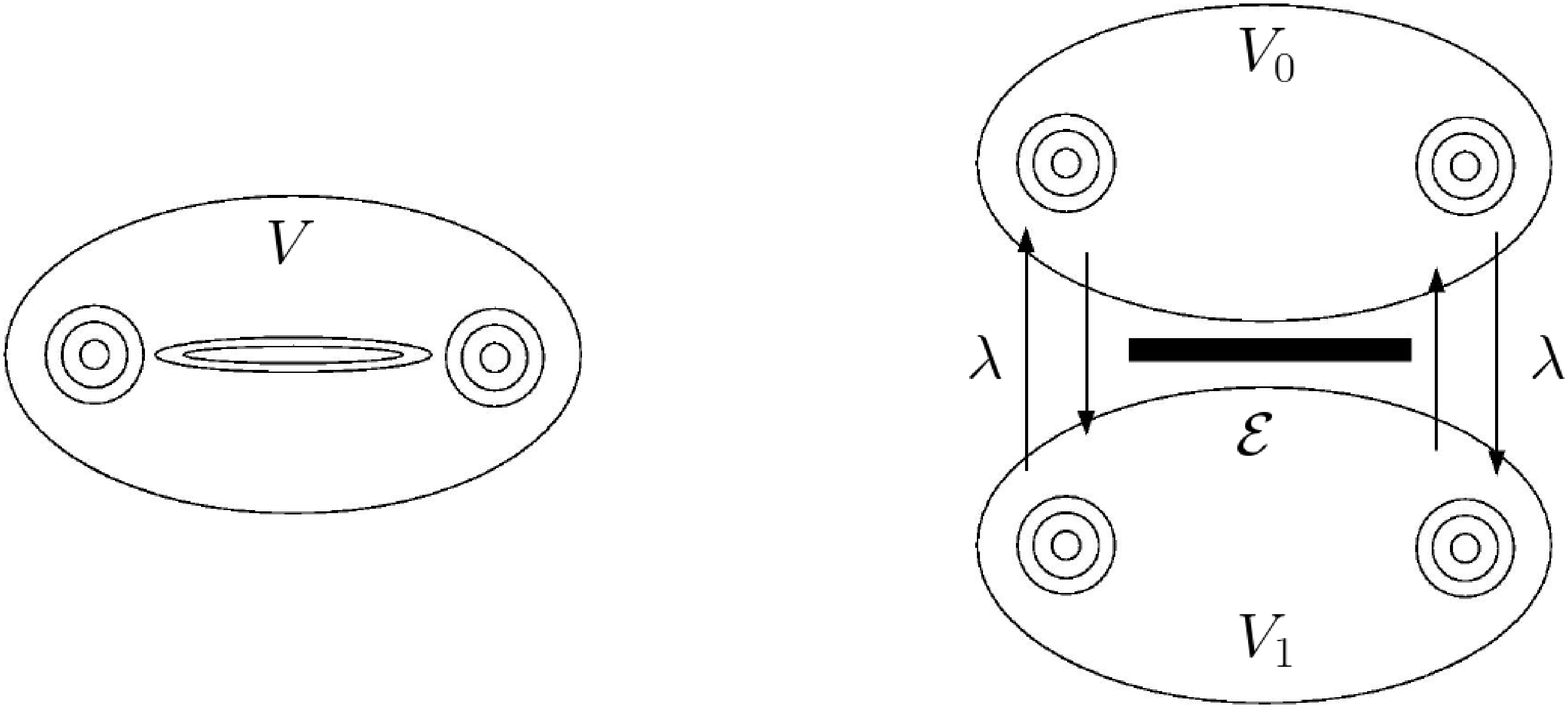}
 \caption{Left: Contour plot of a 2-dimensional potential energy surface demonstrating the bi-channel scenario.  {Right:} In  the bi-channel model the two channels are described by two potential functions $V_i: \D \rightarrow \R$, $i\in \{0,1\}$.  Exchange  between the two channels is allowed only in regions $\T \backslash \e$ at a rate $\lambda > 0$.  }
 \label{fig:bichannel}
 \end{figure}

The bi-channel situation is characterized by the existence of two channels joining an initial and final state on a potential energy surface $V$, separated from each other by a region of high energy, as depicted in the left of Figure~\ref{fig:bichannel}. As explained above, we were not able to 
analyze the original ABF dynamics~\eqref{eq:abf_orig} in this situation, because of some technical difficulties in expressing the probability density flux from one channel to the other.  

We therefore analyze the convergence of the ABF method for a slightly different model, which is schematically represented in right of Figure~\ref{fig:bichannel}.  Each channel is described by a potential energy function $V_i:\D \rightarrow \R$, where $i \in  \{0,1\}$ denotes the channel  index. The stochastic process we consider now is actually a couple $(X_t,I_t)$, where the position vector $X_t$ lives at time $t$ on the potential $V_{I_t}$, $I_t \in \{0,1\}$ being the index of the visited channel at time $t$. The channel index $I_t$ is allowed to switch to $1-I_t$ only if $\xi(X_t)$ lies in some designated regions (typically at the two ends of the two channels).

The dynamics for the ABF dynamics in the bi-channel model is then
\begin{equation}
\label{eq:abf}
\left\{
\begin{aligned}
 dX_t     &= -\nabla(V_{I_t}-A_t \circ \xi)(X_t)dt + \sqrt{2}\ dB_t, \\
 A'_t (x) &= \E\left[\partial_xV_{I_t}(X_t)|\xi(X_t)=x\right], \\
 I_t \in \{0,1\} &\text{ is a jump process with generator}\\
& \text{ $L\varphi(x,y,i) = -\lambda(x) (\varphi(x,y,i)-\varphi(x,y,1-i))$}.
\end{aligned}
\right.
\end{equation}
In terms of the stochastic process $I_t$, switching between the two potentials (namely change of $I_t$ to $1-I_t$) occurs at times 
$$\tau_{n+1} = \text{inf}\left\{~t > \tau_{n} ~\left|~\int_{\tau_{n}}^{t}\lambda(\xi(X_s)) ds > T_n~\right\}\right.,$$ 
where $\tau_0 = 0$ and $T_n$ are i.i.d. random variables drawn from the exponential distribution with parameter 1.   In this way $\lambda(x)$ denotes the rate at which the trajectories jump from potential $V_i$ to potential $V_{1-i}$.  Note that this rate depends only on the position, $x$, in the reaction coordinate.  The bi-channel feature of the model is related to the fact that we assume that the rate $\lambda$ is zero (there is no switching) in some region of the reaction coordinate. Outside of this region, the rate is supposed to be constant, and the potential functions are assumed to be identical (the particles live in the same potential). Let us state this as a formal assumption,
\begin{equation*}
\text{{\bf [H1]}~~~}
\exists~\e \subset \T, \left.
\begin{array}{c}
\text{$\lambda(x) = \lambda \mathbf{1}_{\T \backslash \e}(x)$~\text{and} $\forall x \in \T\backslash \e, V_0(x,\cdot) = V_1(x,\cdot)$}.
\end{array}
\right.
\end{equation*}
The region $\e \subset \T$ in the above hypothesis represents the region where the two channels are separated by high energy barriers, see Figure~\ref{fig:bichannel}. It is assumed to have a Lebesgue measure different from 0 and 1.  

The main qualitative difference between the bi-channel model we study, and the original ABF dynamics~\eqref{eq:abf_orig} is that the switching only depends on the $x$-position and not on the $y$-position. However, the ABF dynamics~\eqref{eq:abf} conserves the main difficulty of the original one, namely the metastability of the dynamics in terms of visited channels for some values of $\xi$. At times $t$ such that $\xi(X_t)\in \e$, $I_t$ cannot switch to $1-I_t$.  In particular, it can be checked that the proof of~\cite{lrs_nonlinearity} applied to~\eqref{eq:abf} in the case $\e = \emptyset$ leads to an estimated rate of convergence limited by $\lambda$ and is thus eventually zero if $\lambda$ goes to zero (see Remark~\ref{rem:lambda_positif} below).  The aim of this work is to study the case $\e \neq \emptyset$ and to obtain an exponential rate of convergence even if $\lambda = 0$ in some region (see Theorem~\ref{th:main} below).  

\begin{remark}
\label{rem:multiple_channels}
On Figure~\ref{fig:bichannel}, we represent the bi-channel case with {\em two} metastable states linked by two different channels. The region where the stochastic process can jump from one channel to another has two connected components. We would like to emphasize that our result also applies to the case where this region is simply connected, namely a situation where two channels start from a metastable state (along the reaction coordinate value) but do not end in another metastable state. This situation is of course less favorable in terms of speed of convergence to equilibrium, which would be reflected in our theoretical result through the parameter $\theta$ (see [H4] below). This is actually typically the situation of the numerical experiments in~\cite{minoukadeh2010} mentioned above since the various compact states obtained do not belong to the same metastable basin.  
\end{remark}

\subsection{A partial differential equation formulation}
\label{sec:PDE}
Let us introduce the time marginal of the process $(X_t,I_t)$ in \eqref{eq:abf}: $$d\mu_t = \sumi \psi(t,x,y,i)\delta_i~dx~dy,$$ where $\delta_i$ is the Dirac measure on the singleton $\{i\}$.  When necessary, we shall denote the $i$-dependency of the density by a subscript $\psi_i$, so that $\psi_{1-i}$ denotes $\psi(t,x,y,1-i)$ for example.  The evolution of the densities are described by a system of non-linear partial differential equations: $\forall i \in \{0,1\}$,
\begin{equation}
 \label{eq:FP}
  \partial_t \psi = \text{div}(\nabla (V_{i} - A_t \circ \xi)\psi) + \Delta \psi -\lambda \circ \xi \, [\psi-\psi_{1-i}] \text{ on $\D$},
\end{equation}
\noindent where, we recall, the last term is zero for $\xi(x,y) \in \e$.   The non-linearity is due to the definition of the mean force estimate, given by 
\begin{equation}
\label{eq:A'} 
\displaystyle A'_t(x) = \frac{\displaystyle \sumi\int_{\R} \partial_x V_i(x,y) \ \psi(t,x,y,i)\ dy}
				{\displaystyle \sumi\int_{\R} \psi(t,x,y,i) \ dy}.   
\end{equation}
\noindent Using hypothesis [H1], it can be checked that if $\psi_{\infty}$ is defined as a probability density proportional to $e^{-(V_i-A_{\infty}\circ \xi)}$ where $A_\infty$ is a given long-time limit for $A_t$, then $\psi_\infty$ is a stationary solution to~\eqref{eq:FP}. Then, by replacing $\psi$ by $\psi_{\infty}$ in \eqref{eq:A'} and comparing with the definition of the mean force \eqref{eq:act_A'}, it is clear that by choosing $A'_{\infty} = A'$, one obtains a stationary solution of the system \eqref{eq:FP}--\eqref{eq:A'} written as:
\begin{equation}
\label{eq:psialpha}
\psi_{\infty}(x,y,i) = \frac{\displaystyle e^{-(V_i(x,y)-A(x))}}{\displaystyle \sumi \int_{\D}e^{-(V_i-A\circ \xi)}~dx~dy}.  
\end{equation}
The associated equilibrium measure for the process $(X_t,I_t)$ writes:
$$d\mu_\infty = \sumi \psi_\infty(x,y,i)\delta_i~dx~dy.$$
We will need further notation for marginal and conditional laws associated to $\mu_t$ and $\mu_\infty$. The next two sections give precise definitions for these measures, and Table~\ref{tab:notation} summarizes the notation.

\begin{table}[t]\centering
\begin{tabular}{|c|c|l|}
\hline
& &  \\[-1ex]
\textbf{Measure}  &\textbf{Definition} &\textbf{Description of measure}    \\[1ex]
\hline 
$d\mu_t $	& $\displaystyle \sumi \psi(t,x,y,i)\delta_i~dx~dy$  
		& Probability measure in domain $\D \times \{0,1\}$ \\[3ex]

$d\mu_t^{\xi,I}$& $\displaystyle \sum_{i=0}^1\psixi(t,x,i) \delta_i \, dx$ 
				& {Marginal measure in reaction coordinate and channel }  \\[3ex]

$d\mu_t^{\xi}$	& $\displaystyle \psixibar(t,x)\, dx$ 
				& Marginal measure in reaction coordinate   \\[3ex]

$d\mu_{t|x,i}$	& $\frac{\displaystyle\psi(t,x,y,i) \ dy}{\displaystyle\psixi(t,x,i)} $ 
			& Measure conditioned to being at $x$ and $i$ \\[3ex]

$d\mu_{t|x}$	& $\frac{\displaystyle\sum_{i=0}^1 \psi(t,x,y,i) \delta_i \ dy}{\displaystyle\psixibar(t,x)}$  
			& Measure conditioned to being at $x$ \\[3ex]

$d\mu^I_{t|x}$	& $ \frac{\displaystyle \sumi \psixi(t,x,i) \delta_i}{\displaystyle\psixibar(t,x)}$ 
			& Marginal in $I$, conditioned to being at $x$ \\[3ex]
\hline
\end{tabular}
\caption{A list of probability measures. }
\label{tab:notation}
\end{table}

\subsubsection{Marginal laws}
\label{sec:marginal_laws}
We are now in the position to define marginal probability densities.  The image of the probability measure $\mu_t$ in $\xi$ and $I$ is denoted by 
$$d\mu_t^{\xi,I}=\sum_{i=0}^1\psixi(t,x,i) \delta_i \, dx ~\text{  where  }~ \psixi(t,x,i) = \int_{\R}\psi(t,x,y,i) \ dy.$$
The evolution of the density $\psixi$ is described by the (non-closed) partial differential equation
\begin{equation}
 \label{eq:FPpsixi}
  \partial _t \psixi = \int_{\R}\partial_x(\partial_x (V_{i} - A_t \circ \xi)\psi) \ dy 
	+ \partial_{xx} \psixi -\lambda(x)(\psixi -\psixi_{1-i}), 
\end{equation}
obtained by integrating \eqref{eq:FP} in $y$.  The associated equilibrium measure is $$d\mu^{\xi,I}_{\infty}=\sum_{i=0}^1\psixi_\infty(x,i) \delta_i \, dx~\text{  where  }~
 \psixi_{\infty}(x,i) = \int_{\R}\psi_{\infty}(x,y,i)~dy.  
$$
\noindent The marginal measure in $\xi$ only is denoted by 
$$d\mu_t^{\xi}=\psixibar(t,x)\, dx ~\text{ where }~\psixibar(t,x) = \sumi\psixi(t,x,i)=\sumi\int_\R \psi(t,x,y,i) \, dy.$$
By summing \eqref{eq:FPpsixi} over $i$ and using the definition~\eqref{eq:A'} of $A_t'$, it is easy to check that $\psixibar$ satisfies a closed, very simple partial differential equation (this is similar to the original ABF dynamics, see~\cite{lrs2007,lrs_nonlinearity}):
\begin{equation}
 \label{eq:FPpsixibar}
  \partial _t \psixibar = \partial_{xx} \psixibar~~\text{on}~\T.  
\end{equation}
Thanks to the adaptive bias, along $\xi$, the barriers have been flattened.  The long-time limit of $\psixibar$ is given by
$$\psixibar_{\infty}(x) = \sumi \int_{\R}\psi_{\infty}(x,y,i)~dy = 1,$$
which corresponds to the uniform probability measure on the torus $\T$.

\subsubsection{Conditional laws}
\label{sec:conditional_laws}
Let us introduce the measure $\mu_{t|x}$ of $(X_t,I_t)$ conditioned to being at a specified point $x$ in the reaction coordinate:
\begin{equation}\label{eq:mutx}
d\mu_{t|x} = \frac{\displaystyle\sum_{i=0}^1 \psi(t,x,y,i) \delta_i \ dy}{\psixibar(t,x)}.
\end{equation}
Its long-time limit is
\begin{equation}\label{eq:muinftytx}
d\mu_{\infty|x} = \frac{\displaystyle\sum_{i=0}^1 \psi_{\infty}(x,y,i)\delta_i  \ dy}{\psixibar_{\infty}(x)}.  \end{equation}
To study the bi-channel model, we will also need to introduce the measures conditioned to being at fixed $\xi$ \textit{and} a particular channel $i \in \{0,1\}$. This measure and its long-time limit are defined respectively by
\begin{equation}
 \label{eq:mu_xi}
d\mu_{t|x,i} = \frac{\psi(t,x,y,i) \ dy}{\psixi(t,x,i)} 
		\ \ \ \text{and} \ \ \   
d\mu_{\infty|x,i} = \frac{\psi_{\infty}(x,y,i) \ dy}{\psixi_{\infty}(x,i)}.  
\end{equation}
Finally, we will also need the marginals in $I$ of the probability measures $\mu_{t|x}$ and $\mu_{\infty|x}$. These Bernoulli probability measures $\mu^I_{t|x}$ and $\mu^I_{\infty|x}$ represent the proportion of the marginal distribution $\psixibar(t,\cdot)$ and $\psixibar$ in each channel.  They are formally defined as
\begin{equation}
 \label{eq:mu_i}
d\mu^I_{t|x} =  \frac{\displaystyle \sumi \psixi(t,x,i) \delta_i}{\psixibar(t,x)}
		\ \ \ \text{and} \ \ \   
d\mu^I_{\infty|x} =  \frac{\displaystyle \sumi \psixi_\infty(x,i) \delta_i}{\psixibar_\infty(x)}.
\end{equation}

\subsection{Entropy and Fisher information}
\label{sec:entropy_defs}
In this section we recall some well-known results for defining relative entropy between two probability measures, which can be seen as a measure of the `distance' between those. A general introduction to this topic can be found in~\cite{logsob2000,villanitopics} and applications to study the longtime behavior of Fokker-Planck type equations are presented in~\cite{arnold2001}.  
\begin{definition} \textbf{(Entropy).}
For any two probability measures $\mu$ and $\nu$ such that $\mu$ is absolutely continuous with respect to $\nu$ (denoted as $\mu \ll \nu$), the relative entropy is defined as\end{definition}
$$H(\mu|\nu) = \int \ln \left( \frac{d\mu}{d\nu}\right) d\mu.  $$
The positivity of the relative entropy can be shown using the inequality $\ln(1/x) \geq 1-x$.  Furthermore, $H(\mu|\nu) = 0$ if and only if $\mu = \nu$.  
\begin{definition} \textbf{(Csiszar-Kullback inequality). }
For measures $\mu$ and $\nu$ which have densities with respect to the Lebesgue measure, the following holds
\begin{equation*}\label{eq:csiszar}\Vert\mu-\nu\Vert_{L^1} \leq \sqrt{2H(\mu|\nu)}.  \end{equation*}
\end{definition}
This allows us to control the $L^1$-norm of the difference of two probability measures by their relative entropy.  
\begin{definition}\textbf{(Fisher information).}
 For any probability measure $\mu$ absolutely continuous with respect to $\nu$, the Fisher information is given by \end{definition}
$$F(\mu|\nu) = \int \left| \nabla \ln \left( \frac{d\mu}{d\nu} \right)\right|^2 d \mu.  $$
\begin{definition}\textbf{(Logarithmic Sobolev inequality).}
 The probability measure $\nu$ is said to satisfy a logarithmic Sobolev inequality with constant $\rho > 0$ (in short: LSI$(\rho)$) if for all probability measures~$\mu$ such that $\mu \ll \nu$,
$$H(\mu|\nu) \leq \frac{1}{2\rho}F(\mu|\nu).$$
\end{definition}
\noindent Such an inequality holds for Gaussian measures~\cite{gross1975} for example, and more generally~\cite{bakry1984} for any measure $\nu$ with density proportional to $e^{-V}$, where $V$ is $\alpha$-convex (in which case the LSI constant $\rho$ is equal to $\alpha$). Besides, there exists a perturbation result~\cite{holleystroock}: if $\tilde{\nu}$ is a probability measure such that $d \tilde{\nu}/d \nu=e^{U}$, where $\nu$ satisfies a LSI($\rho$) and $U$ is a bounded function, then $\tilde{\nu}$ satisfies a LSI with constant $\tilde{\rho}=\rho \, e^{-\text{osc}(U)}$, where 
$\text{osc}(U)=\sup(U)-\inf(U)$. Thus, a very large class of probability measures satisfy a LSI. An important feature of the LSI constant is that it degenerates to zero in the case where the probability measure is multimodal. For example, if $d\nu= Z^{-1} \exp(-\beta W(x)) \, dx$ and $W(x)=x^4/4-x^2/2$ is the double-well potential in dimension 1, then the LSI constant scales as $\exp(-\beta \Delta W)$ where $\Delta W=W(0)-W(1)>0$ is the height of the barrier. Such inequalities thus hold under rather loose assumptions, but the constant $\rho$ is very small for a multimodal measure. For example, in typical situations encountered in molecular dynamics, the LSI constant for the measure $\nu$ defined in \eqref{eq:langevin} is extremely small, which is related to the fact that it is difficult to sample directly this measure. 

To analyze the convergence of the ABF dynamics for the bi-channel model, we make use of LSIs for the conditional  measures $\mu_{\infty|x,i}$ (see assumption [H3] below), which are the equilibrium canonical measures, conditioned to being at a fixed value of $\xi$ {\em and in a channel}. In~\cite{lrs_nonlinearity}, log-Sobolev inequalities for the equilibrium canonical measures conditioned only to being at a fixed value of $\xi$ were considered, but in the bi-channel case, those are typically very small due to the presence of high energy barriers `orthogonal' to the isocontours of the reaction coordinate (see again Figure~\ref{fig:bichannel}).

Let us now define the Wasserstein distance between two probability measures.  
\begin{definition}\textbf{(Wasserstein distance).}
The Wasserstein distance with linear cost between probability measures $\mu$ and $\nu$ is defined as 

$$\displaystyle W(\mu,\nu) = \inf_{\pi \in \Pi(\mu,\nu)} \int_{\D \times \D}
			|y-y'| \ \pi(dy,dy'),$$
where $\Pi(\mu,\nu)$ denotes the set of coupling probability measures on $\D \times \D$, with marginals $\mu$ and $\nu$.  
\end{definition}
\begin{definition}\textbf{(Talagrand inequality).}
 The probability measure $\nu$ is said to satisfy a Talagrand inequality with constant $\rho > 0$ (or $T(\rho)$) if for all probability measures $\mu$ such that $\mu \ll \nu$, 
\begin{equation} 
\label{eq:talagrand}
W(\mu,\nu) \leq \sqrt{\frac{2}{\rho}H(\mu|\nu)}. 
\end{equation}
\end{definition}
Logarithmic Sobolev inequalities and Talagrand inequalities are related (see~\cite{otto2000}):
\begin{lemma}
\label{lem:T}
 If $\nu$ satisfies $LSI(\rho)$, then $\nu$ satisfies $T(\rho)$.  
\end{lemma}

Below, we present entropies that prove useful in obtaining convergence results of the bi-channel ABF model.  
In this paper, we are primarily interested in the convergence to a stationary state of the Fokker-Planck equation \eqref{eq:FP}--\eqref{eq:A'}, and thus of the associated partial differential equations \eqref{eq:FPpsixi} and \eqref{eq:FPpsixibar}.  Relative entropies will therefore be defined for some probability measures with respect to their long-time limits. 

The \textit{total entropy} will be denoted by
$$E(t) = H(\mu_t|\mu_{\infty}) = \sumi \int_{\D}\ln\left(\frac{\psi}{\psi_{\infty}}\right)\psi \ dx \ dy.  $$
The so-called \textit{macroscopic entropy} is defined as
$$E_M(t) 	= H(\mu_t^{\xi}|\mu^{\xi}_{\infty}) 
		= \int_{\T}\ln\left(\frac{\psixibar}{\psixibar_{\infty}}\right)\psixibar \ dx , $$
the \textit{local entropy} at a fixed value $x$ in the reaction coordinate by
$$e_m(t,x) = H(\mu_{t|x}|\mu_{\infty|x})=	\sumi \int_{\R} \ln\left(\frac{\psi}{\psixibar} \Bigg/
				\frac{\psi_{\infty}}{\psixibar_{\infty}}\right) \frac{\psi}{\psixibar} \ dy,$$
and the \textit{microscopic entropy} by
$$E_m(t) = \int_{\T} e_m(t,x) \ \psixibar(t,x) \ dx.$$
With the above definitions, it is easy to show that
\begin{equation}
\label{eq:sumE}
E(t) = E_M(t) + E_m(t).
\end{equation}
In order to treat the bi-channel case, we define a \textit{channel-local entropy}, defined by 
$$ e_{\text{cl}}(t,x,i) = 	H(\mu_{t|x,i}|\mu_{\infty|x,i}) = 
				\int_{\R} \ln\left(\frac{\psi}{\psixi} \Bigg/
				\frac{\psi_{\infty}}{\psixi_{\infty}}\right) \frac{\psi}{\psixi} \ dy. $$

Two hypothesis that will be essential to the results presented below are: a so-called `bounded coupling' assumption on the cross derivative $\partial_{x,y}V_i$ (see~\cite{lrs_nonlinearity,lelievre2010}), and an assumption on the logarithmic Sobolev inequality constants for the probability measures $\mu_{\infty|x,i}$. The hypotheses read
\begin{equation*}
\text{{\bf [H2]}~~~}
\left\{
\begin{array}{c}
\text{$\forall i \in \{0,1\}$, $V_i$ and $\xi$ are sufficiently differentiable functions such that~$\exists~C, M > 0$} \\
\text{$\left\| \partial_{x,y}V_i \right\|_{L^\infty(\T \times \R)} \leq M<\infty$ and  
$\displaystyle \left\|\frac{\int_{\R}\partial_xV_i e^{-V_i}~dy}
			{\int_{\R}e^{-V_i}~dy}\right\|_{L^{\infty}(\T)} \leq C < \infty$,  }
\end{array}
\right.
\end{equation*}
\begin{equation*}\label{eq:hyp_LSI_orig}
\text{{\bf [H3]}~~~}
\left\{
\begin{array}{c}
\text{$\forall i \in \{0,1\}$, $V_i$ and $\xi$ are such that $\exists \rho >0$, $\forall x \in
  \T$, $\forall i \in \{0,1\}$,}\\
\text{the conditional measure $\mu_{\infty|x,i}$ satisfies LSI($\rho$).}
\end{array}
\right.
\end{equation*}
\noindent The hypothesis [H3] gives us that $\forall x \in \T,  \forall i \in \{0,1\}$
\begin{eqnarray*}
 H( \mu_{t|x,i} | \mu_{\infty|x,i} ) & \leq & \frac{1}{2\rho} F(\mu_{t|x,i} | \mu_{\infty|x,i} ),  
\end{eqnarray*}
\noindent or equivalently
\begin{equation}
\int_{\R} \ln\left( \frac{\psi}{\psixi } \Bigg/ \frac{\psi_\infty }{
\psixi_\infty }\right) \frac{\psi }{\psixi  } \  dy \\
\leq  \frac{1}{2\rho} \int_{\R} \left| \partial_y \ln \left( \frac{\psi}{\psi_{\infty}}\right) \right|^2
                             \frac{\psi}{\psixi  } \ dy  \label{lsi_i}.  \\
\end{equation}
Finally, an entropy that appears in later calculations, is that of the Bernoulli measure $\mu^I_{t|x}$ with respect to its long-time limit
$$ e_c(t,x) = H(\mu^I_{t|x}|\mu^I_{\infty|x}) 
	= \sumi \ln \left( \frac{\psixi}{\psixibar}\Bigg/\frac{\psixiinf}{\psixibar_{\infty}} \right) 
		\frac{\psixi}{\psixibar}.  $$
The so-called \textit{channel entropy} is then given by 
\begin{eqnarray}
 E_c(t) &=& \int_{\T} e_c(t,x) \psixibar \ dx  \nonumber \\
	&=&  \sumi \int_{\T}   \ln \left( \frac{\psixi}{\psixibar}\Bigg/\frac{\psixiinf}{\psixibar_{\infty}} \right) \psixi.  \label{eq:Ec_def} 
\end{eqnarray}

\subsection{The free energy as a bias in each channel}
\label{sec:good_bias}
As well as the hypothesis [H3] on the conditional measures $\mu_{\infty|x,i}$, an assumption will also be necessary to ensure that the free energy  is a `good bias' {\em in each channel}. More precisely, once the bias has converged, the marginals along $\xi$ in each channel must converge sufficiently quickly to their long-time limit. Roughly speaking, the channels should not be too `asymmetrical'. The aim of this section is to state this more formally, see assumption [H4] below.

 Consider that the system is already nearly at equilibrium, in the sense that 
\begin{enumerate}
 \item[i)] $A_t'=A_{\infty}'=A'$,
 \item[ii)] $\forall i \in \{0,1\}$, $\displaystyle \int_{\R}\partial_xV_i \ d\mu_{t|x,i}  = \int_{\R}\partial_xV_i \ d\mu_{\infty|x,i}$.
\end{enumerate}
Then the marginal density $\psixi$ can be shown to satisfy (see~\eqref{eq:FPpsixi} and~\eqref{eq:fplem2} below)
\begin{equation}
\label{eq:psixi_op}
\partial_t\psixi = -\L_i \psixi
	= \partial_x \left(\psixi_\infty \partial_x (\psixi / \psixi_\infty) \right)
		-\lambda(x)(\psixi -\psixi_{1-i}).  
\end{equation}
It can be shown that the operator $\L = (\L_0,\L_1)$ is symmetric and positive definite with respect to the inner product 
$\langle f, g \rangle = \displaystyle \sumi\int_{\T} f_i(x) g_i(x) \ \frac{1}{\psixi_{\infty}(x,i)}\ dx$ and has a spectral gap $\theta > 0$ (see Lemma~\ref{lem:theta} below).  
In other words, $\L$ is such that for all functions $f:\T\times\{0,1\}\rightarrow \R$, with $f_i \in H^1\left(\frac{1}{\psixi_{\infty}(x,i)}dx\right)$, $f \neq 0$ and 
$\sumi\int_{\T} {f_i(x)}~ dx =0$ , we have
\begin{equation}
\label{eq:poincare}
 \langle f,f \rangle 
	\leq 
 \frac{1}{\theta}\langle f, \L f \rangle.  
\end{equation}

As will be apparent in the proof, it will be necessary for the spectral gap $\theta$ to be sufficiently large.  
\begin{equation*}\label{eq:marg_min}
\text{{\bf [H4]}~~~}
\left\{
\begin{array}{c}
\text{$\displaystyle \theta > \theta_{\rm min}$ with $\displaystyle \theta_{\rm min}=\frac{8(C+M\rho^{-1/2})^2\tM}{c}$}\\[5pt]
\text{where $0 < c,\tM < \infty$ are such that $\displaystyle \inf_{x,i} \psixi_{\infty} = c$ and $\sup_{x} \psixibar(0,x) = \tM $}.
\end{array}
\right.
\end{equation*}
\noindent We recall that the constants $C$, $M$ and $\rho$ have been introduced in assumptions [H2] and [H3] above. The fact that $c>0$ is not restrictive since $\forall i\in \{0,1\}$, $V_i$ is a continuous function and $\psixi_{\infty} \propto \int_{\R}e^{-V_i}\ dy > 0$ is continuous and defined on the compact space $\T$.  Similarly, since $\psixibar$ satisfies the heat equation~\eqref{eq:FPpsixibar}, the assumption $\tM < \infty$ is not restrictive. If, for example, the initial condition has a Dirac mass marginal in $\xi$, for any positive time $t_0>0$, $\psixibar(t_0,\cdot)$ is a bounded function, and one has simply to consider the dynamics on $[t_0, \infty)$.

\begin{remark}
In the hypothesis [H2], the assumption on the cross derivative $\partial_{x,y} V_i$ could in fact be replaced by $\|\partial_xV_i\|_{L^{\infty}(\T \times \R)} \leq M < \infty$, in which case, in [H4], the minimum value for $\theta$ would be  changed to $\theta_{\rm min}=\frac{8M^2\tM}{c}$ and the Talagrand inequalities in Lemmas \ref{lem:At_A_Em} and \ref{lem:diff_xi_bound} would be replaced by Csiszar-Kullback inequalities.  
\end{remark}

\subsection{Main result}
\label{sec:results}
We are now in position to present the main result of the paper.  
\begin{theorem}
\label{th:main}
Assume hypotheses [H1]-[H4]. There exists a smooth function $\Lambda: (\theta_{\rm min}, \infty) \rightarrow (0,\rho)$ which is increasing and such that:
\begin{equation*}
\Lambda(\rho+2\theta_{\rm min}) = \frac{\rho}{2}
~\text{ and }~
  \Lambda(\theta)  \rightarrow
\left\{
\begin{array}{r l}
 0      &\text{as}~~\theta \rightarrow \theta_{\rm min}\\
  \rho &\text{as}~~\theta \rightarrow \infty
\end{array}
\right.
\end{equation*}
for which we can prove the following convergence result: for any $\varepsilon \in (0, \Lambda(\theta))$, there exists a constant $K>0$ such that, $\forall t>0$, 
\begin{equation}\label{eq:CVEm}
 E_m(t) \leq K \text{{\em exp}} \left(-2 \, \min\{(\Lambda(\theta)-\varepsilon),4\pi^2\} \, t\right).
\end{equation}
This implies that the total entropy $E$ and thus $\|\psi(t,\cdot)-\psi_{\infty}\|^2_{L^1}$ converge exponentially fast to zero with the same rate.  
Furthermore, the biasing force $A_t'$ converges to the mean force  $A'$ in the following sense: $\forall t \ge 0$,
$$\int_{\T}|A_t'(x) - A'(x)|^2 \psixibar(t,x) \ dx \leq 2(C+M\rho^{-1/2})^2 \,E_m(t).$$
As a consequence, for any positive time $t_0>0$ and $\varepsilon \in (0, \Lambda(\theta))$, there exists a constant~$\bar{K}$ such that $\forall t \ge t_0$,
$$\int_{\T}|A_t'(x) - A'(x)|^2 \ dx \leq \bar{K} \text{{\em exp}} \left(-2 \, \min\{(\Lambda(\theta) - \varepsilon),4\pi^2\} \, t\right).$$
\end{theorem}
The term $4 \pi^2$ corresponds to the exponential rate of convergence of $\psixibar$ to $\psixibar_\infty$ (see~\eqref{eq:FPpsixibar}), and is clearly not the bottleneck. There are actually various ways to make this rate as small as needed (see Remark~11 in~\cite{lrs_nonlinearity}).

Thus, this result essentially shows that the ABF method converges at a rate which is limited by the multimodality of the equilibrium canonical measures conditioned to being at fixed value of $\xi$ {\em and in a channel} (quantified by the constant $\rho$), if the free energy is a bias which enables a fast exploration of each channel (this is quantified by the constant $\theta$, which should be sufficiently large for $\Lambda(\theta)$ to be indeed close to $\rho$). Thus, the convergence may be fast even if switching between the two channels is impossible for some values of $\xi$. If the spectral gap $\theta$ is sufficiently large, we thus recover a similar expression for the rate of convergence of the ABF method as the one derived in~\cite{lrs_nonlinearity}, with $\rho$ being now the LSI constant of the canonical measures $\mu_{|x,i}$.

We would like to emphasize that our arguments hold under the following assumption of existence of regular solutions: We assume that we are given a process $(X_t,I_t)$ and a function $A_t'$ which satisfy~\eqref{eq:abf}, and such that $X_t$, conditionally on $I_t=i$ has a smooth density $\psi(t,x,y,i)$. We suppose that this density is sufficiently regular so that the entropy estimates below are valid. We refer for example to~\cite{arnold2001} for an appropriate functional framework in which such entropy estimates hold.

\section{Proof of main result}
\label{sec:proof}
In order to prove the exponential decay of the microscopic entropy in Theorem \ref{th:main}, we use the fact that the time evolution of the microscopic entropy can be expressed as a combination of the evolution of the total and the marginal entropies, from \eqref{eq:sumE}
$$ \frac{dE_m}{dt} = \frac{dE}{dt} - \frac{dE_M}{dt}.  $$
In order to obtain results for the microscopic entropy, we begin by treating the time evolutions of the total entropy and the macroscopic entropy separately.  


\subsection{Preliminary computations on the total entropy $E$}
In order to study the evolution of the total entropy, some auxiliary results will be needed and are given in the lemmas below.  First, it will be useful to write the Fokker-Planck equation associated to $\psi$ in a different form.  

\begin{lemma}
\label{lem:FP2}
The Fokker-Planck equation~\eqref{eq:FP} for $\psi$ can be rewritten as
$$ \partial_t\psi 
	= \div( \psi_\infty \nabla (\psi / \psi_\infty)) 
		+ \partial_x( (A' - A_t') \psi)
		-\lambda(x)(\psi -\psi_{1-i}). $$
\end{lemma}
\proof
By developing the right hand side, we obtain
\begin{eqnarray*}
 \partial_t\psi 
	&=&  \div \left(\nabla \psi - \frac{\psi}{\psi_{\infty}} \nabla\psi_{\infty} \right) 
		+ \partial_x((A'-A_t')\psi) 
		-\lambda(x)(\psi -\psi_{1-i}) \\
	&=& \div \left(\nabla \psi + \nabla(V_i-A\circ \xi)\psi\right) 
		+ \partial_x(A'\psi) - \partial_x(A_t'\psi)  
		-\lambda(x)(\psi -\psi_{1-i})\\
	&=& \div \left(\nabla(V_i-A_t\circ \xi)\psi \right) + \Delta \psi 
		-\lambda(x)(\psi -\psi_{1-i}),
\end{eqnarray*}
which is indeed the Fokker-Planck equation \eqref{eq:FP}.  \endproof

\noindent The above may now be used to estimate the time evolution of the total entropy.  

\begin{lemma}
\label{lem:global}
 The total entropy $E$ satisfies
\end{lemma}
\begin{equation}
\label{eq:glob_ineq}
 \frac{dE}{dt} 
	\leq  -\sumi \int_{\T \times \R} 
		\left| \nabla\ln\left( \frac{\psi}{\psi_{\infty}}\right) \right|^2 \psi
		- \sumi \int_{\T \times \R} (A'-A'_t) \partial_x 
                    \left[\ln\left( \frac{\psi}{\psi_{\infty}}\right)\right]\psi. 
\end{equation}
\proof
First, by definition of the total entropy, we have
\begin{eqnarray}
 \label{eq:ode_E}
\displaystyle
 \frac{dE}{dt} & = & \frac{d}{dt} \sumi \int_{\T \times \R} \ln\left( \frac{\psi}{\psi_{\infty}}\right) \psi \nonumber \\
& = & \sumi \int_{\T \times \R} \ln \left( \frac{\psi}{\psi_{\infty}} \right)\partial_t\psi,  
\end{eqnarray}

\noindent using the fact that $\psi$ is a probability density.  Next, we use Lemma \ref{lem:FP2} to obtain

\begin{eqnarray}
 \frac{dE}{dt} 
	&=& \sumi \int_{\T \times \R}\text{div}\left(\psi_{\infty} \nabla \left(\frac{\psi}{\psi_{\infty}}\right)\right)
		\ln\left( \frac{\psi}{\psi_{\infty}}\right) 
		+ \sumi \int_{\T \times \R} \partial_x((A'-A'_t) \psi ) 
		\ln\left( \frac{\psi}{\psi_{\infty}}\right) \nonumber \\
	& & - \sumi \int_{\T \times \R} \lambda(x) \ \ln \left( \frac{\psi}{\psi_{\infty}} 
		\right) (\psi-\psi_{1-i}) \label{eq:lam_negative} \\ 
	&\leq& -\sumi \int_{\T \times \R} 
		\left| \nabla\ln\left( \frac{\psi}{\psi_{\infty}}\right) \right|^2 \psi
		- \sum_{i=0}^1\int_{\T \times \R} (A'-A'_t) \partial_x 
		\left[\ln\left( \frac{\psi}{\psi_{\infty}}\right)\right]\psi,\nonumber 
\end{eqnarray}
where the last line is a result of integration by parts and the fact that \eqref{eq:lam_negative} is non-positive, which is proved below in Lemma~\ref{lem:g2}.  \endproof
\begin{lemma}
 \label{lem:g2}
 For $\psi$ satisfying \eqref{eq:FP} and $\psi_{\infty}$ its long-time limit, the following holds \end{lemma}
$$ -\sumi \int_{\T \times \R} \lambda(x) \ \ln \left( \frac{\psi}{\psi_{\infty}} \right) (\psi-\psi_{1-i}) \leq 0. $$
\proof
 First, recall that $\lambda(x) = 0$, $\forall x \in \e$.  We consider therefore the left hand side of the inequality at fixed $(x,y) \in (\T\backslash \e) \times \R$.  At fixed $(x,y)$, we consider renormalized (Bernoulli) probabilities denoted by $\tilde{\psi}_i~=~\psi_i/(\psi_0+\psi_1)$, so that $\tpsiz + \tpsio = 1$. 
\begin{align*}
-\sumi &\ln\left(\frac{\psi}{\psi_{\infty}}\right)(\psi-\psi_{1-i})\\
	&= (-\psi_0 + \psi_1 ) \left[ \ln \left(\frac{\psi_0}{\psi_{\infty,0}}\right)
		- \ln \left(\frac{\psi_1}{\psi_{\infty,1}} \right) \right] \\
	&= (\psi_0+\psi_1)(- \tpsiz + \tpsio )
		 \left[ \ln \left(\frac{\tpsiz}{\tpsizinf}\right)- \ln \left(\frac{\tpsio}{\tpsioinf} \right) \right] \\
	&= (\psi_0+\psi_1)\left[- \ln \left(\frac{\tpsiz}{\tpsizinf}\right)\tpsiz
		- \ln \left(\frac{\tpsio}{\tpsioinf} \right)\tpsio \right. 
	 	+ \left.  \ln \left(\frac{\tpsiz}{\tpsizinf}\right)\tpsio
		+  \ln \left(\frac{\tpsio}{\tpsioinf} \right) \tpsiz   \right] \\      
	&= (\psi_0+\psi_1)\left[- 2\ln \left(\frac{\tpsiz}{\tpsizinf}\right)\tpsiz
		- 2\ln \left(\frac{\tpsio}{\tpsioinf} \right)\tpsio \right. 
	 	+ \left. (\tpsiz+\tpsio) \ln \left(\frac{\tpsiz}{\tpsizinf}\right) 
		+ (\tpsiz+\tpsio)  \ln \left(\frac{\tpsio}{\tpsioinf} \right)    \right] \\      
	&= (\psi_0+\psi_1) \left[ - 2\sum_{i=0}^{1}\ln \left(\frac{\tpsizi}{\tpsizinfi}\right)\tpsizi
		+ \ln \left(\frac{\tpsiz}{\tpsizinf}\right)
		+ \ln \left(\frac{\tpsio}{\tpsioinf} \right) \right].  
\end{align*}
From hypothesis [H1], $\forall (x,y)\in (\T\backslash \e) \times \R$, $\tpsioinf = \tpsizinf = \frac{1}{2}$, which allows the above to be written as
\begin{eqnarray}
 -\sumi \ln\left(\frac{\psi}{\psi_{\infty}}\right)(\psi-\psi_{1-i}) 
	&=&  (\psi_0+\psi_1)\left[ -2\sum_{i=0}^{1}\ln \left(\frac{\tpsizi}{\tpsizinfi}\right)\tpsizi    
		 + 2\sum_{i=0}^1  \ln \left(\frac{\tpsizi}{\tpsizinfi}\right)\tpsizinfi  \right] \nonumber \\
	&=& - 2(\psi_0+\psi_1) \left[ \sum_{i=0}^{1}\ln \left(\frac{\tpsizi}{\tpsizinfi}\right)\tpsizi
		+ \sum_{i=0}^1  \ln \left(\frac{\tpsizinfi}{\tpsizi}\right)\tpsizinfi \right] \label{eq:inv_entropies} \\
	& \leq &  0.  \nonumber 
\end{eqnarray}
The last line is due to the fact that the two terms between brackets are non-negative, since they are relative entropies.  
\endproof

\begin{remark}\label{rem:lambda_positif}
Let us consider the case where $\e = \emptyset$ and thus $\forall x \in \T,~ \lambda(x) = \lambda >0$, which implies (see [H1]) $V_0 = V_1$ and therefore $\psi_{\infty,0}=\psi_{\infty,1}$ everywhere on $\T\times\R$. In this case, it follows from~\eqref{eq:inv_entropies} that
\begin{align*}
 -\sumi\int_{\T\times\R} & \lambda(x) 
                \text{{\em ln}} \left(\frac{\psi}{\psi_{\infty}}\right)(\psi-\psi_{1-i}) \\
  &\leq - 2\lambda\int_{\T\times\R}(\psi_0+\psi_1) \sum_{i=0}^{1}\text{{\em ln}} \left(\frac{\tpsizi}{\tpsizinfi}\right)\tpsizi \\
  &\leq - 2\lambda\int_{\T\times\R}\sum_{i=0}^{1}\text{{\em ln}} \left(\frac{\psi}{\psi_{\infty}}\right)\psi
	+ 2\lambda\int_{\T\times\R}\text{{\em ln}} \left(\frac{\psi_0+\psi_1}{\psi_{\infty,0}+\psi_{\infty,1}}\right)(\psi_0+\psi_1).  
\end{align*}
\def\psixy{\psi^{x,y}}
Furthermore, since the marginal $\psixy := \psi_0+\psi_1$ satisfies in this specific case
\begin{equation}
\partial_t\psixy 	= \div( \psixy_\infty \nabla (\psixy / \psixy_\infty))
		+ \partial_x( (A' - A_t') \psixy),
\end{equation}
one can show, using the results of~\cite{lrs_nonlinearity}, that $\exists C > 0$, $\forall t \geq 0$,
$$ -\sumi\int_{\T\times\R}   \lambda(x) \text{{\em ln}} \left(\frac{\psi}{\psi_{\infty}}\right)(\psi-\psi_{1-i}) 
    \leq 
  - 2\lambda E +2\lambda C \text{{\em e}}^{-2 \, \min\{\rho,4\pi^2\} \, t} . $$
Using this, and the fact that $-E \leq -E_m$, one can show that $E_m$ converges to zero exponentially fast with rate
$$ 2\min\{\rho, 4\pi^2,\lambda\}.  $$
The convergence rate thus depends on $\lambda > 0$, the rate at which switching occurs between the two channels.  This is comparable to the original result obtained for the ABF algorithm, see~\cite{lrs_nonlinearity}.  The aim of what follows is to obtain a result in the case where $\e \neq \emptyset$, namely when $\lambda = 0$ in some region.  
\end{remark}

\noindent In order to estimate the last term on the right hand side of~\eqref{eq:glob_ineq}, it will be helpful to express the difference of the biasing force and the mean force in terms of densities

\begin{lemma}
\label{lem:At_A}
The difference between the biasing force $A_t'$ and the mean force $A'$ can be expressed in the following way
\end{lemma} 
$$A_t'-A'= \sumi \int_{\R} \partial_x \ln\left(\frac{\psi}{\psi_{\infty}}\right)\frac{\psi}{\psixibar} \ dy
	- \partial_x \ln\left( \frac{\psixibar}{\psixibar_{\infty}}\right).$$
\proof We develop the expression on the right hand side and use the fact that $\psixibar_{\infty}\equiv 1$

\begin{align*}
\sumi \int_{\R} \partial_x \ln\left(\frac{\psi}{\psi_{\infty}}\right)\frac{\psi}{\psixibar} \ dy
	&- \partial_x \ln\left( \frac{\psixibar}{\psixibar_{\infty}}\right)  \\
	&= 	\sumi \int_{\R} \partial_x \ln(\psi) \frac{\psi}{\psixibar} \ dy
		- \sumi \int_{\R} \partial_x \ln(\psi_{\infty}) \frac{\psi}{\psixibar} \ dy 
		- \partial_x \ln\left(\psixibar \right) \\
	&= 	\sumi \int_{\R} \frac{\partial_x \psi}{\psixibar} \ dy
		+ \sumi \int_{\R} \partial_x (V_i -A\circ\xi) \frac{\psi}{\psixibar} \ dy 
		- \frac{\partial_x \psixibar} {\psixibar} \\
	&=	\sumi \int_{\R} \partial_x V_i  \frac{\psi}{\psixibar} \ dy
		- \sumi \int_{\R} A' \frac{\psi}{\psixibar} \ dy  \\
	&=	A_t' - A',
\end{align*}
The last line is a result of the definition of $A_t'$ in \eqref{eq:A'} and the fact that $A$ is a function of $x$ only.  \endproof 
Another useful estimate for the difference between $A_t'$ and $A'$ is given in the following lemma.  
\begin{lemma}
\label{lem:At_A_Em}
The difference of the biasing force and the mean force can be bounded by the microscopic entropy as:
 $$\int_{\T} |A_t' - A'|^2 \psixibar \ dx \leq  2\tR^2E_m(t),$$
where $$\tR = \left(C+M\rho^{-1/2}\right).$$  
\end{lemma}
\proof 
We begin by showing that $|A_t'(x)-A'(x)| \leq M\sqrt{2e_m(t,x)/\rho}$.  By definition, we have
\begin{eqnarray*}
\displaystyle
A_t'(x)-A'(x) 
   &=& 
        \frac{\displaystyle \sumi \int_{\R}\partial_xV_i \psi \ dy}{\displaystyle\sumi \int_{\R}\psi \ dy} 
            - \frac{\displaystyle \sumi \int_{\R}\partial_xV_i \psi_{\infty} \ dy}{\displaystyle\sumi \int_{\R}\psi_{\infty} \ dy} \\
   &=& \sumi \int_{\R} \left( \left( \partial_xV_i \frac{\psixi}{\psixibar} \right) \frac{\psi}{\psixi}
            -\left(\partial_xV_i  \frac{\psixi_{\infty}}{\psixibar_{\infty}}\right)
			 \frac{\psi_{\infty}}{\psixi_{\infty}}  \right)   \ dy \\
   &=& \sumi \int_{\R\times\R} \left(  \partial_xV_i(x,y)\frac{\psixi}{\psixibar} 
       -\partial_xV_i(x,y')  \frac{\psixi_{\infty}}{\psixibar_{\infty}}\right)
			 \pi(dy,dy') \\
\end{eqnarray*}
where $\pi(dy,dy')$ is any coupling measure on $\R\times\R$ with marginals $\mu_{t|x,i}$ and $\mu_{\infty|x,i}$.  Next, using Taylor's expansion on $\partial_xV_i(x,y)$, we have
\begin{eqnarray*}
\displaystyle
A_t'(x)-A'(x) 
   &=& \sumi \int_{\R\times\R} 
	\left(  \left(\partial_xV_i(x,y') + \partial_{x,y}V_i(x,\eta(y,y'))(y-y')\right)\frac{\psixi}{\psixibar} 
       -\partial_xV_i(x,y')  \frac{\psixi_{\infty}}{\psixibar_{\infty}}\right) \pi(dy,dy') \\
   &=& \sumi \int_{\R\times\R} 
	\left(  \partial_xV_i(x,y')\left(\frac{\psixi}{\psixibar}-\frac{\psixi_{\infty}}{\psixibar_{\infty}}\right)
	+ \partial_{x,y}V_i(x,\eta(y,y'))(y-y')\frac{\psixi}{\psixibar}\right)  \pi(dy,dy') \\
   &=&    \sumi \left(\frac{\psixi}{\psixibar}-\frac{\psixi_{\infty}}{\psixibar_{\infty}}\right)
	\int_{\R} \partial_xV_i(x,y')\frac{\psi_{\infty}}{\psixi_{\infty}}\ dy' 
	+ \sumi\frac{\psixi}{\psixibar}  \int_{\R\times\R} \partial_{x,y}V_i(x,\eta(y,y'))(y-y') \pi(dy,dy') \\
\end{eqnarray*}
where $\eta(y,y') \in [y,y']$.  Recall from [H2] that $\exists C, M > 0$ such that $\displaystyle \left\|\int_{\R}\partial_xV_i\ d\mu_{\infty|x,i}\right\|_{L^{\infty}(\T)} \leq C$ and $\|\partial_{x,y}V_i\|_{L^{\infty}} \leq M$.  Furthermore, with the use of the Csiszar-Kullback inequality, we have

\begin{eqnarray*}
\displaystyle
|A_t'(x)-A'(x)| 
   &\leq&    C \sumi \left|\frac{\psixi}{\psixibar}-\frac{\psixi_{\infty}}{\psixibar_{\infty}}\right|
	+ M  \sumi \frac{\psixi}{\psixibar}\int_{\R\times\R}|y-y'| \pi(dy,dy') \\
   &\leq&   C \sqrt{2 H(\mu^I_{t|x},\mu^I_{\infty|x})}
	+ M  \sumi \frac{\psixi}{\psixibar}\int_{\R\times\R}|y-y'| \pi(dy,dy').  
\end{eqnarray*}
Next, by the Talagrand inequality \eqref{eq:talagrand} and the concavity of the function $x\mapsto \sqrt{x}$, 
\begin{eqnarray*}
|A_t'(x)-A'(x)| 
   &\leq&    C \sqrt{2 H(\mu^I_{t|x},\mu^I_{\infty|x})}
	+ M \sumi \frac{\psixi}{\psixibar}\sqrt{\frac{2}{\rho}H(\mu_{t|x,i},\mu_{\infty|x,i})} \\
   &\leq&    C \sqrt{2H(\mu_{t|x},\mu_{\infty|x})}
	+ M \sqrt{\frac{2}{\rho}\sumi H(\mu_{t|x,i},\mu_{\infty|x,i})\frac{\psixi}{\psixibar}} \\
   &\leq&    C \sqrt{2H(\mu_{t|x},\mu_{\infty|x})}
	+ M \sqrt{\frac{2}{\rho}H(\mu_{t|x},\mu_{\infty|x})} \\
\end{eqnarray*}
Therefore
$$ |A_t'(x)-A'(x)|^2  \leq   2\left(C+M\rho^{-1/2}\right)^2 e_m(t,x). $$
The result follows immediately, since $E_m(t) = \int_{\T}e_m(t,x) \ \psixibar dx$.   \endproof

\noindent We are now equipped with the right tools to control $dE/dt$, and are left to handle the evolution of the macroscopic entropy.

\subsection{Controlling $E_M$}
Due to the free diffusion equation~\eqref{eq:FPpsixibar} satisfied by $\psixibar$, the macroscopic entropy $E_M$ is easily controlled.  

\begin{lemma}
\label{lem:macro}
The macroscopic entropy satisfies \end{lemma}
\begin{equation}
 \label{eq:ode_EM} 
\frac{dE_M}{dt}	= - \int_{\T} \left| \partial_x \ln 
			\left( \frac{\psixibar}{\psixibar_{\infty}} \right) \right|^2 \psixibar
		= - F(\psixibar|\psixibar_{\infty}).  
\end{equation}
\proof Using \eqref{eq:FPpsixibar} and integration by parts
\begin{eqnarray*}
 \frac{dE_M}{dt} &=& \frac{d}{dt}\int_{\T} \ln \left( \frac{\psixibar} {\psixibar_{\infty}} \right)\psixibar \\
                 &=& \int_{\T} \ln \left( \frac{\psixibar}{\psixibar_{\infty}} \right) \partial_{xx} \psixibar \\
		 &=& - \int_{\T} \left| \partial_x \ln \left( \frac{\psixibar}{\psixibar_{\infty}} \right) \right|^2 \psixibar.    
\end{eqnarray*}

\begin{lemma}
 \label{lem:macro_fisher}
The Fisher information for the marginal density $\psixibar$ decays exponentially fast with rate $r = 8\pi^2$ 
$$F(\psixibar(t,\cdot)|\psixibar_{\infty}) \leq F_0~\text{{\em exp}}(-8\pi^2 t),  $$
where $F_0 =  F(\psixibar(0,\cdot)|\psixibar_{\infty})$.  
\end{lemma}
\proof See Lemma 12 of reference~\cite{lrs_nonlinearity}. 

 \endproof

\noindent In light of the lemmas presented above, hypotheses [H1]-[H4] may now be used to control the evolution of the microscopic entropy.  

\subsection{Controlling $E_m$}
The aim of this section is to obtain an estimation on the evolution of th microscopic entropy $E_m$, see~\eqref{eq:dEm_dt} below.  
We first begin by using the fact that the channel-local conditional measures $\mu_{\infty|x,i}$ satisfy LSI($\rho$) to bound the microscopic entropy.  
\begin{lemma}
\label{lem:E_m}
 Under the hypothesis [H3], the microscopic entropy $E_m$ satisfies \end{lemma}
\begin{equation*}
 E_m  \leq \frac{1}{2\rho} \sum_{i=0}^1\int_{\T \times \R}  \left| \partial_y \ln \left(
                    \frac{\psi}{\psi_{\infty}}  \right) \right|^2   \psi  \ dx \ dy
            + E_c,
\end{equation*}
{\em where $E_c$ is defined in \eqref{eq:Ec_def}.  } \\
\proof 
We fist consider the local entropy $e_m(t,x)$, which can be decomposed into the sum of the channel-local entropy and the entropy of the Bernoulli measures.    
\begin{eqnarray*}
 e_m(t,x) & = &  H( \mu_{t|x} | \mu_{\infty|x} )  \\
          & = & \sumi \int_{\R} \ln \left(\frac{\psi}{\psixibar} \Big/ \frac{\psi_{\infty}}{\psixibar_{\infty}} \right)
			 \frac{\psi}{\psixibar} \ dy   \\
          & = & \sumi \left[ \int_{\R}  \ln \left(
                    \frac{\psi}{\psixi}   \Big/ \frac{\psi_{\infty}}{\psixiinf}  \right)
                            \frac{\psi}{\psixibar}  \ dy 
                +   \int_{\R} \ln \left( \frac{\psixi }{\psixibar}   \Big/ \frac{\psixiinf}{\psixibar_{\infty}}  
                    \right)      \frac{ \psi}{\psixibar} \ dy \right]   \\
          & \leq & \sumi \left[ \frac{1}{2\rho} \int_{\R}  \left| \partial_y \ln \left(
                    \frac{\psi}{\psi_{\infty}}  \right) \right|^2
                            \frac{\psi}{\psixibar} \ dy  \right]
                +  \sumi \left[   \ln \left( \frac{\psixi}{\psixibar}   \Big/ \frac{\psixiinf}{\psixibar_{\infty}}  
                    \right)  \frac{\psixi}{\psixibar} \right] \\
          & = & \frac{1}{2\rho} \sumi  \int_{\R}  \left| \partial_y \ln \left(
                    \frac{\psi}{\psi_{\infty}}  \right) \right|^2
                            \frac{\psi}{\psixibar} \ dy 
                +  e_c(t,x), \\
\end{eqnarray*}
where the inequality is a direct result of [H3].  The microscopic entropy is then 
\begin{eqnarray*}
\label{eq:Em}
\displaystyle
 E_m(t) &=&	\int_{\T} e_m(t,x) \psixibar\  dx  \\
	&\leq& \frac{1}{2\rho} \sum_{i=0}^1\int_{\T \times \R}  \left| \partial_y \ln \left(
		\frac{\psi}{\psi_{\infty}}  \right) \right|^2   \psi  \ dx \ dy
		+ E_c,  
\end{eqnarray*} 
as required.  \endproof

\noindent The time evolution of $E_m$ may now be expressed using results of Lemmas \ref{lem:global} and \ref{lem:macro}.  
\begin{eqnarray*}
 \frac{dE_m}{dt} 
	&=&  	\frac{dE}{dt}- \frac{dE_M}{dt} \\
	&\leq&  	-\sum_{i=0}^1\int_{\T \times \R}\left| \nabla\ln\left( \frac{\psi}{\psi_{\infty}}\right) \right|^2 \psi 
		- \sumi \int_{\T \times \R} (A'-A'_t) \partial_x \ln\left( \frac{\psi}{\psi_{\infty}}\right)\psi  
		+ \int_{\T} \left| \partial_x \ln \left( \frac{\psixibar}{\psixibar_{\infty}} \right) \right|^2 \psixibar.  
\end{eqnarray*}
We may now apply Lemma \ref{lem:At_A} and use integration by parts  
\begin{eqnarray}
 \frac{dE_m}{dt} 
	&=& 	-\sumi \int_{\T \times \R}\left| \partial_y \ln\left( \frac{\psi}{\psi_{\infty}}\right) \right|^2 \psi 
		-\sumi \int_{\T \times \R}\left| \partial_x \ln\left( \frac{\psi}{\psi_{\infty}}\right) \right|^2 \psi \nonumber \\
	& & 	+\sumi \int_{\T \times \R} \left(\sumi \int_{\R} \partial_x
			\ln\left(\frac{\psi}{\psi_{\infty}}\right)\frac{\psi}{\psixibar} \ dy
		- \partial_x \ln\left( \frac{\psixibar}{\psixibar_{\infty}}\right) \right) \partial_x 
			\ln\left( \frac{\psi}{\psi_{\infty}}\right)\psi \nonumber \\ 
	& & 	+ \int_{\T} \left| \partial_x \ln \left( \frac{\psixibar}{\psixibar_{\infty}} \right) \right|^2 \psixibar  \nonumber \\
	&=& 	-\sumi \int_{\T \times \R}\left| \partial_y \ln\left( \frac{\psi}{\psi_{\infty}}\right) \right|^2 \psi \nonumber \\
	& &	-\sumi \int_{\T \times \R}\left| \partial_x \ln\left( \frac{\psi}{\psi_{\infty}}\right) \right|^2 \psi 
		+\int_{\T}\left(\sumi \int_{\R}\partial_x 
		\ln\left(\frac{\psi}{\psi_{\infty}}\right)\psi\ dy \right)^2\frac{1}{\psixibar} \label{eq:cauchy-schwartz} \\
	& &	-\sumi \int_{\T \times \R}\partial_x \ln\left( \frac{\psixibar}{\psixibar_{\infty}}\right)  \partial_x 
		\ln\left( \frac{\psi}{\psi_{\infty}}\right)\psi 
		+\int_{\T} \left| \partial_x \ln \left( \frac{\psixibar}{\psixibar_{\infty}} \right) \right|^2 \psixibar.   \nonumber
\end{eqnarray}
Notice that \eqref{eq:cauchy-schwartz} is non-positive by the Cauchy-Schwartz inequality.  We therefore have 
\begin{eqnarray}
\frac{dE_m}{dt} 
	&\leq&	-\sumi \int_{\T \times \R}\left| \partial_y \ln\left( \frac{\psi}{\psi_{\infty}}\right) \right|^2 \psi \nonumber \\
	& &	- \sumi \int_{\T \times \R}\partial_x \ln\left( \frac{\psixibar}{\psixibar_{\infty}}\right)  \partial_x 
		\ln\left( \frac{\psi}{\psi_{\infty}}\right)\psi  
		+ \int_{\T} \left| \partial_x \ln \left( \frac{\psixibar}{\psixibar_{\infty}}\right) \right|^2 \psixibar \nonumber \\
	&=& 	-\sumi \int_{\T\times\R}\left| \partial_y \ln\left(\frac{\psi}{\psi_{\infty}}\right) \right|^2 \psi   
		-\int_{\T }\left[\partial_x \ln\left( \frac{\psixibar}{\psixibar_{\infty}}\right)\psixibar \right] (A_t'-A') \label{eq:at_a} \\
	&\leq&	-2\rho E_m  
		+2\rho E_c
		+ \sqrt{\int_{\T}\left|A_t'-A'\right|^2\psixibar}\sqrt{\int_{\T}
		\left| \partial_x \ln\left(\frac{\psixibar}{\psixibar_{\infty}}\right)\right|^2\psixibar}.  \nonumber 
\end{eqnarray}
Line \eqref{eq:at_a} is a result of Lemma \ref{lem:At_A} and the last inequality is due to Lemma \ref{lem:E_m} and a further application of the Cauchy-Schwartz inequality.  Now, using Lemmas \ref{lem:At_A_Em} and \ref{lem:macro_fisher}, we obtain
\begin{eqnarray*}
 \frac{dE_m}{dt} &\leq& -2\rho E_m  
			+2\rho E_c
                        + \tR\sqrt{2E_m}\sqrt{F_0\text{e}^{-8\pi^2t}},.  
\end{eqnarray*}
where we recall $F_0 = F(\psixibar(0,\cdot)|\psixibar_{\infty})$.  Finally, using Young's inequality: $\forall \varepsilon > 0$, $\forall a,b \in \R$, $ab < \varepsilon a^2 + \frac{1}{4\varepsilon}b^2$, we obtain
\begin{eqnarray}
\label{eq:dEm_dt}
 \frac{dE_m}{dt} &\leq& -2\left(\rho - \tR^2\varepsilon \right) E_m  
			+2\rho E_c
			+ \frac{1}{4\varepsilon}F_0\text{e}^{-8\pi^2t},  
\end{eqnarray}
where $\varepsilon > 0$ will be chosen optimally later in the proof.  We are left to control the channel entropy term $E_c$ in order to conclude.  

\subsection{Controlling $E_c$}
\label{sec:g1}
The aim of this section is to obtain a control on the evolution of $E_c$ defined by~\eqref{eq:E_c}, and more precisely an upper bound of $E_c$, denoted $P$, which is the weighted $\chi^2$-distance between $\psixi$ and $\psixi_{\infty}$, see~\eqref{eq:E_c} and~\eqref{eq:firstP} below.  

Recall that $E_c = \int_{\T}H(\mu^I_{t|x}|\mu^I_{\infty|x})\psixibar \ dx$, where the integrand is the relative entropy of a Bernoulli measure.  Poincar\'e and logarithmic Sobolev inequalities have been studied for Bernoulli measures~\cite{logsob2000}.  In order to obtain an exponentially decaying relative entropy, however, a suitable semi-group and its infinitesimal generator is needed for the measure.  In the case of the bi-channel model (in particular [H1]), we face problems due to the region $\e \subset \T$, where no exchange is permitted between the two channels: in this region, the speed at which the measure $\mu^I_{t|x}$ reaches equilibrium cannot directly be controlled.  To circumvent this issue, we consider the spectral gap of an adequate operator and resort to the Poincar\'e inequality.  

By the definition of $E_c$ and using the inequality $\forall x > 0$, $x\ln(x) \leq x(x-1)$ and the fact that 
$\sumi \int_\T \psixi(t,x,i) \, dx = 1$, we obtain
\begin{eqnarray}
E_c
	&=&   \sumi \int_{\T}\ln \left(\frac{\psixi}{\psixibar}\Bigg/
                           \frac{\psixi_{\infty}}{\psixibar_{\infty}}\right)\psixi \nonumber \\
	&=&    \sumi \int_{\T}\ln \left(\psixi/\psixi_{\infty}\right)\psixi  - E_M   \nonumber \\
	&\leq& 	\sumi \int_{\T}\left(\frac{\psixi}{\psixi_{\infty}}-1\right)^2\psixi_{\infty} \ dx 
		- E_M.  
\end{eqnarray}

\noindent We therefore have
\begin{equation} E_c \leq P, \label{eq:E_c}\end{equation}
where
\begin{equation}
\label{eq:P}
P = \sumi \int_{\T}\left(\frac{\psixi}{\psixi_{\infty}}-1\right)^2\psixi_{\infty} \ dx.  
\end{equation}
In order to proceed and consider the time derivative of $P$, we will need some further results to express the evolution of the marginal density $\psixi$.  The idea is to compare the evolution of $\psixi$ with the dynamics of this density if $A_t'$ and 
$\int_{\T}\partial_xV d\mu_{t|x,i}$ were already at equilibrium (see Section \ref{sec:good_bias}).  

\begin{lemma}
\label{lem:diff_xi}
The Fokker-Planck equation \eqref{eq:FPpsixi} for $\psixi$ can be rewritten as
\begin{equation}
\label{eq:fplem2}
\begin{array}{l l l}
\partial_t\psixi 
	&=& \displaystyle \partial_x \left(\psixi_\infty \partial_x \left(\frac{\psixi}{\psixi_\infty}\right) \right)
		+ \partial_x\left(\left( \frac{\int_{\R}\partial_xV_i \psi}{\psixi} 
		- \frac{\int_{\R}\partial_x V_i \psi_{\infty}}{\psixi_{\infty}}\right) \psixi\right) \\
	& & + \displaystyle \partial_x( (A' - A_t') \psixi)
		-\lambda(x)(\psixi -\psixi_{1-i}) 
\end{array}
\end{equation}
\end{lemma}
\proof First we will show that
\begin{equation}
\label{eq:fplem1}
\partial_t\psixi = \partial_x \left(\int_{\R} \psi_\infty \partial_x \left(\frac{\psi}{ \psi_\infty}\right)\ dy\right) 
	+ \partial_x( (A' - A_t') \psixi)
	-\lambda(x)(\psixi -\psixi_{1-i}).  
\end{equation}
By developing the right hand side, we have
\begin{eqnarray*}
 \partial_t\psixi 
     &=& \partial_x \left(\int_{\R} \partial_x \psi -\frac{\psi}{\psi_{\infty}}\partial_x{\psi_{\infty}} \right) 
         +  \partial_x \int_{\R} (A'-A_t')\psi   -\lambda(x)(\psixi -\psixi_{1-i})\\
     &=& \partial_x \left(\int_{\R}\partial_x(V_i-A\circ \xi)\psi \right) +\partial_{xx} \psixi
         +  \partial_x \int_{\R} (A'-A_t')\psi   -\lambda(x)(\psixi -\psixi_{1-i})\\
     &=& \int_{\R} \partial_x \left(\partial_x(V_i-A_t\circ \xi)\psi \right)
              +\partial_{xx} \psixi -\lambda(x)(\psixi -\psixi_{1-i}) \\
\end{eqnarray*}
which is indeed the Fokker-Planck equation \eqref{eq:FPpsixi} associated to $\psixi$.  

Next, we show that $\forall x \in \T$, $\forall i \in \{0,1\}$
\begin{equation}
\label{eq:lemterm}
\int_{\R} \psi_\infty \partial_x  \left(\frac{\psi}{ \psi_\infty}\right)\ dy 
	=
  \psixi_\infty \partial_x \left(\frac{\psixi}{\psixi_\infty}\right)
  + \left( \frac{\int_{\R}\partial_xV_i \psi}{\psixi} - \frac{\int_{\R}\partial_xV_i \psi_{\infty}}{\psixi_{\infty}}\right) \psixi.   
\end{equation}
To prove the above, notice that
\begin{eqnarray*}
 \int_{\R} \psi_\infty \partial_x  \left(\frac{\psi}{ \psi_\infty}\right)\ dy 
	&-&  \psixi_\infty \partial_x \left(\frac{\psixi}{\psixi_\infty}\right) \\
	&=&	\int_{\R}\partial_x \psi \ dy - \int_{\R}\frac{\psi}{\psi_{\infty}}\partial_x\psi_{\infty}
		- \partial_x \psixi + \frac{\psixi}{\psixi_{\infty}}\partial_x \psixi_{\infty} \\
	&=& 	\int_{\R}(\partial_x(V_i-A\circ \xi)) \psi \ dy 
		- \frac{\displaystyle\int_{\R}\partial_x(V_i-A\circ \xi)e^{-(V_i-A\circ \xi)} dy} 
		       {\displaystyle \int_{\R}e^{-(V_i-A \circ \xi)}\ dy} \psixi \\
	&=& \int_{\R}(\partial_xV_i) \psi\ dy -A'(x)\psixi 
		- \frac{\displaystyle \int_{\R}\partial_x V_i e^{-(V_i-A\circ \xi)}\ dy}
		 {\displaystyle \int_{\R}e^{-(V_i-A \circ \xi)}\ dy} \psixi + A'(x)\psixi
\end{eqnarray*}
Finally, by using the fact that the free energy $A$ is independent of $y$, we obtain
\begin{equation}
  \int_{\R} \psi_\infty \partial_x (\psi / \psi_\infty)\ dy 
		-  \psixi_\infty \partial_x (\psixi / \psixi_\infty)
	= \left( \frac{\int_{\R}\partial_xV_i \psi}{\psixi} 
	  - \frac{\int_{\R}\partial_xV_i\psi_{\infty}}{\psixi_{\infty}}\right) \psixi,   
\end{equation}
as required.  The final result \eqref{eq:fplem2} is obtained by substituting \eqref{eq:lemterm} into \eqref{eq:fplem1}.  \endproof

\noindent Notice that \eqref{eq:fplem2} is comparable to~\eqref{eq:psixi_op}, only with additional terms due to the fact that $A_t$ and 
$\int_{\T}\partial_xV d\mu_{t|x,i}$ have not yet converged.  The difference of the biasing force and mean force, $A_t'-A'$, was already estimated in Lemmas \ref{lem:At_A} and \ref{lem:At_A_Em}.  We are therefore left to control the remaining term.  

\begin{lemma}
\label{lem:diff_xi_bound}
$\forall x \in \T$, $\forall i \in \{0,1\}$,
$$ \left| \frac{\int_{\R}\partial_xV_i\psi}{\psixi}-\frac{\int_{\R}\partial_xV_i\psi_{\infty}}{\psixi_{\infty}}\right|      \leq M\sqrt{\frac{2}{\rho}H\left( \mu_{t|x,i}\left| \mu_{\infty|x,i}\right.\right)}.  $$
As a consequence, 
\begin{equation}
\label{eq:diff_xi}
\displaystyle
	\sumi \int_{\T} 
	\left| \frac{\int_{\R}\partial_xV_i\psi}{\psixi} 
	- \frac{\int_{\R}\partial_xV_i\psi_{\infty}}{\psixi_{\infty}}\right|^2  \psixi \ dx
	\leq \frac{2M^2}{\rho}E_m.  
\end{equation}
\end{lemma}
\proof Let $\Pi(\mu_{t|x,i},\mu_{\infty|x,i})$ be the set of coupling measures on $\R \times \R$ with marginals $\mu_{t|x,i}$ and $\mu_{\infty|x,i}$ respectively and let $\pi \in \Pi$ .  Then
\begin{eqnarray*}
 \left|\int_{\R}\left(\partial_xV_i\frac{\psi}{\psixi}-\partial_xV_i\frac{\psi_{\infty}}{\psixi_{\infty}}\right)\right| 
	&=& 
		\left|\int_{\R\times \R}\left(\partial_xV_i(x,y)-\partial_xV_i(x,y')\right)\ \pi(dy,dy')
			\right| \\
	&\leq& \left\|\partial_{x,y}V_i\right\|_{L^{\infty}}\int_{\R\times \R}|y-y'|\ \pi(dy,dy') \\
	&\leq& M\sqrt{\frac{2}{\rho}H(\mu_{t|x,i}| \mu_{\infty|x,i})} \\
\end{eqnarray*}
\noindent where we have used Lemma \ref{lem:T} since $\mu_{\infty|x,i}$ satisfies LSI$(\rho)$.  Equation \eqref{eq:diff_xi} follows immediately from the fact that $\sumi \int_{\T}H(\mu_{t|x,i}| \mu_{\infty|x,i}) \psixi \ dx \leq E_m$.  \endproof

One result that is now needed to derive estimates on the evolution of $P$ is the existence of a spectral gap of the operator describing the dynamics once $A'_t$ and $\displaystyle \int_{\R}\partial_xV_id\mu_{t|x,i}$ have converged (see~\eqref{eq:psixi_op}).  We now justify the existence of such a spectral gap.

In order to do so, let us define the vector spaces
$$\V_l = \left\{ v:\T \times \{0,1\} \to \R~\left|~\forall i \in \{0,1\},
	~\frac{v_i}{\psixi_{\infty}(x,i)}  \in L^2\left(\T, {\psixi_{\infty}(x,i)}\,dx\right)\right. 
,~\sumi\int_{\T}v_i(x) \, dx= l
	\right\}$$
and 
$$\W_l =     \left\{ w \in \V_l~\left|~\forall i \in \{0,1\},
        ~\frac{w_i}{\psixi_{\infty}(x,i)} \in H^1\left(\T, {\psixi_{\infty}(x,i)}~dx\right),~\sumi\int_{\T}w_i(x) \, dx= l \right.
    \right\}. $$
A function $\phi$ in $\V_l$ (or in $\W_l$) will also be considered as a vector valued function as $\phi : \left\{\begin{array}{l} \T \to \R^2 \\ x \mapsto (\phi_0(x),\phi_1(x)) \end{array} \right.$.  Notice that $\phi \in \W_1$ if and only if $f:= \phi-\psixi_{\infty} \in \W_0$.  
\begin{lemma}
\label{lem:theta} 
Recall the operator $\L = (\L_0,\L_1)$, with $\L_i$ defined as in~\eqref{eq:psixi_op},
$$\L_i \phi     =    - \left[\partial_x\left(\psixi_{\infty,i}\partial_x\left(\phi_i/\psixi_{\infty,i}\right)\right)
            - \lambda(x)(\phi_i - \phi_{1-i})\right]. $$
Then
\begin{enumerate}
 \item[i)] The operator $\L$ is symmetric and positive definite with respect to the inner product
$$\displaystyle \langle f, g \rangle = \sumi\int_{\T} f_i(x) g_i(x) \ \frac{1}{\psixi_{\infty}(x,i)}\ dx.$$ 
 \item[ii)] $\L$ has a spectral gap $\theta > 0$ in the sense that 
\begin{equation}
\label{eq:spec_gap}
  \inf_{f \in \W_0, f \neq 0} \frac{\langle f, \L f\rangle}{\langle f,f \rangle} = \theta > 0.
\end{equation}
\end{enumerate}
\end{lemma}

\proof
{\em i)} To show symmetry of the operator $\L$, consider functions $\varphi, \phi~\in~\W_0$. 
Now, using the fact that $\forall x~\in~\T~\backslash~\e$, $\psixi_{\infty,0}(x) = \psixi_{\infty,1}(x)$,
\begin{align}
\sumi &\int_{\T}\varphi_i \L_i \phi \frac{1}{\psixi_{\infty,i}}\ dx \nonumber \\
  &  =    -\sumi \int_{\T}\varphi_i
            \partial_x \left(\psixi_{\infty,i} \partial_x\left(\frac{\phi_i}{\psixi_{\infty,i}}\right)\right)
              \frac{1}{\psixi_{\infty,i}}\ dx
        + \sumi \int_{\T}\lambda(x)\varphi_i(\phi_i - \phi_{1-i})\frac{1}{\psixi_{\infty,i}}\ dx \nonumber \\
    &=    \sumi \int_{\T}\partial_x\left(\frac{\varphi_i}{\psixi_{\infty,i}}\right)
            \partial_x\left(\frac{\phi_i}{\psixi_{\infty,i}}\right) \psixi_{\infty,i}\ dx
        + \int_{\T}\lambda(x)(\varphi_0-\varphi_{1})(\phi_0 - \phi_{1})\frac{1}{\psixi_{\infty,0}}\ dx \label{eq:posdef} \\
   &=    -\sumi \int_{\T} \phi_i  \partial_x \left(\psixi_{\infty,i} \partial_x\left(\frac{\varphi_i}{\phi_{\infty,i}}\right)\right)
            \frac{1}{\psixi_{\infty,i}}\ dx
        + \int_{\T}\lambda(x)(\phi_0 - \phi_{1})(\varphi_0-\varphi_{1})\frac{1}{\psixi_{\infty,0}}\ dx \nonumber \\
    &=    \sumi \int_{\T} \phi_i  \L_i \varphi\frac{1}{\psixi_{\infty,i}}\ dx. \nonumber
\end{align}
From \eqref{eq:posdef}, we conclude positive definiteness of $\L$
\begin{equation}\label{eq:posdef2}
\sumi \int_{\T}\phi_i \L_i \phi \frac{1}{\psixi_{\infty,i}}\ dx
    =    \sumi \int_{\T}\left|\partial_x\left(\frac{\phi_i}{\psixi_{\infty,i}}\right)\right|^2 \psixi_{\infty,i}\ dx
        + \int_{\T}\lambda(x)(\phi_0 - \phi_{1})^2\frac{1}{\psixi_{\infty,0}}\ dx > 0. 
\end{equation}
Notice that the above is strictly positive for any $\phi \in \W_0$, since  $\langle\phi,\L\phi\rangle = 0$ if and only if $\phi_0~=~\phi_1~=~0$.

{\em ii)} In fact, one can check that $\exists \kappa > 0$, such that $\forall \phi \in \W_0$, $\phi \neq 0$,  
$$\sumi \int_{\T} \phi_i \L_i\phi \frac{1}{\psixi_{\infty}}~dx 
       \geq 
  \kappa \sumi \int_{\T} \left(|\phi_i|^2 + \left| \nabla \phi_i \right|^2 \right)\frac{1}{\psixi_{\infty}}~dx. $$
Therefore, by the Lax-Milgram theorem, $\L^{-1}$ is well defined from $\V_0$ to $\W_0$ and thus compact from $\V_0$ to $\V_0$.  
From the symmetry and positive definiteness of $\L$, and the fact that its inverse is a compact operator from $\V_0$ to $\V_0$, it has a strictly positive and discrete spectrum.  There exists a set of eigenvectors $(v_n)_{n\geq 1}$, orthonormal with respect to the inner product $\langle\cdot,\cdot\rangle$, forming a basis of $\V_0$ and $\W_0$, and associated to an increasing sequence of eigenvalues $(\sigma_n)_{n \geq 1}$, such that $\displaystyle \lim_{n \rightarrow \infty}\sigma_n = \infty$.   
In particular, there exists a spectral gap: $\theta=\sigma_1 > 0$.  \endproof

\begin{remark}
 In the case where a function $\phi \in \W_1$ satisfies $\partial_t \phi= \L \phi$, a consequence of Lemma~\ref{lem:theta} is
\begin{equation}
 \label{eq:spec_gap_decay}
\forall t \geq 0,~~\|\phi(t,\cdot)-\psixi_{\infty}\|^2 \leq Ke^{-2\theta t},
\end{equation}
where $\displaystyle\| \cdot \|^2 = \langle\cdot,\cdot\rangle$ and $\displaystyle K=\sum_{n\geq 1} \langle\phi(0,\cdot)-\psixi_{\infty},v_n\rangle^2$.  This is easily obtained by noticing that $\phi-\psixi_{\infty}\in \W_0$ and therefore can be expressed in terms of the orthonormal eigenvectors $(v_n)_{n \geq 1}$  
$$\displaystyle \phi(t,\cdot)-\psixi_{\infty} = \sum_{n\geq 1} \langle \phi(0,\cdot)-\psixi_{\infty},v_n\rangle v_ne^{-\sigma_nt}.  $$ 
The result~\eqref{eq:spec_gap_decay} follows immediately since $\forall n \geq 2$, $\sigma_n \geq \sigma_1 = \theta$.  

Notice in particular that the Fokker-Planck equation~\eqref{eq:fplem2} satisfied by $\psixi$ can be written as
$$\partial_t\psixi(t,x,i) = \L_i\psixi 
    + \partial_x\left(\left( \frac{\int_{\R}\partial_xV_i \psi}{\psixi} 
		- \frac{\int_{\R}\partial_x V_i \psi_{\infty}}{\psixi_{\infty}}\right) \psixi\right) 
    +  \partial_x( (A' - A_t') \psixi.  $$
We will show that the last two terms can be controlled by $E_m$, see~\eqref{eq:firstP}.   
\end{remark}

\noindent With these tools at hand, let us consider the time evolution of the functional $P$ defined in \eqref{eq:P}
\begin{eqnarray*}
 \frac{1}{2}\frac{d P}{dt}
	&=& \sumi \int_{\T}\left(\ffi - 1 \right)\partial_t \left(\ffi\ \right)  \psixiinf \ dx \\
	&=& \sumi \int_{\T}\psixi\partial_t \left(\ffi\right) \  \ dx - \sumi \int_{\T}\partial_t \psixi \ dx\\
	&=& \sumi \int_{\T}\psixi\partial_t \psixi \frac{1}{\psixi_{\infty}} \   dx.  
\end{eqnarray*}
Using equation \eqref{eq:fplem2} of Lemma \ref{lem:diff_xi}, we get
\begin{eqnarray*}
\frac{1}{2}\frac{dP}{dt} 
	&=& \sumi \int_{\T}\psixi\partial_t \psixi \frac{1}{\psixi_{\infty}} \   dx \\
	&=& 	- \sumi\int_{\T} \left|\partial_x \left(\ffi\right)\right|^2\psixiinf \ dx
		- \sumi \int_{\T} \lambda(x)\psixi_i(\psixi_i-\psixi_{1-i})\frac{1}{\psixi_{\infty,i}} \ dx \\
	&&	- \sumi\int_{\T} \partial_x \left(\ffi\right) (A'-A_t')\psixi \ dx
		- \sumi\int_{\T} \partial_x\left(\ffi\right) \left(\frac{\int_{\R}\partial_xV\psi} {\psixi} 
				- \frac{\int_{\R}\partial_xV\psi_{\infty}} {\psixi_{\infty}} \right)\psixi \ dx.  
\end{eqnarray*}
Notice that, by developing the sum and using the fact that $\psixi_{\infty,0}=\psixi_{\infty,1}$ for $\lambda(x) \neq 0$, the second term may be replaced by $\displaystyle \int_{\T} \lambda(x)\left|\psixi_0-\psixi_1\right|^2\frac{1}{\psixi_{\infty,0}} \ dx$.  Finally by using Young's inequality on the last two terms, we obtain for a parameter $\alpha > 0$ to be chosen later on, 
\begin{eqnarray*}
\frac{1}{2}\frac{dP}{dt} 
	&=& 	- \sumi\int_{\T} \left|\partial_x \left(\ffi\right)\right|^2\psixiinf \ dx
		- \int_{\T} \lambda(x)\left|\psixi_0-\psixi_1\right|^2\frac{1}{\psixi_{\infty,0}} \ dx \\
	&&	+ \frac{1}{4\alpha}\sumi\int_{\T} \left|\partial_x \left(\ffi\right)\right|^2 \psixi \ dx
		+ \alpha \int_{\T} |A'-A_t'|^2\psixibar \ dx \\
	&&	+ \frac{1}{4\alpha}\sumi\int_{\T} \left|\partial_x\left(\ffi\right)\right|^2 \psixi \ dx 
		+ \alpha \sumi\int_{\T} \left|\frac{\int_{\R}\partial_xV\psi} {\psixi} 
				- \frac{\int_{\R}\partial_xV\psi_{\infty}} {\psixi_{\infty}} \right|^2\psixi \ dx.  
\end{eqnarray*}
\noindent Next, by Lemmas \ref{lem:At_A_Em} and \ref{lem:diff_xi_bound},
\begin{eqnarray*}
\frac{1}{2}\frac{dP}{dt} 
	&=& 	- \sumi\int_{\T} \left|\partial_x \left(\ffi\right)\right|^2\psixiinf \ dx
		- \int_{\T} \lambda(x)\left|\psixi_0-\psixi_1\right|^2\frac{1}{\psixiinf} \ dx \\
	&&	+ \frac{1}{2\alpha}\sumi\int_{\T} \left|\partial_x \left(\ffi\right)\right|^2 \psixi \ dx 
		+ 2\alpha \frac{M^2}{\rho}E_m  + 2\alpha\tR^2 E_m.  
\end{eqnarray*}
Notice that, in the third term, $\psixi \leq \psixibar \leq \tM$ for $\tM = \left\|\psixibar(0,\cdot)\right\|_{L^{\infty}}$ and $1 \leq \psixi_{\infty}/c$ for $\displaystyle c~=~\min_{x,i}~\psixiinf$.  This gives
\begin{eqnarray*}
\frac{1}{2}\frac{dP}{dt} 
	&\leq& 	- \sumi\int_{\T} \left|\partial_x \left(\ffi\right)\right|^2\psixiinf \ dx
		- \int_{\T} \lambda(x)\left|\psixi_0-\psixi_1\right|^2\frac{1}{\psixiinf} \ dx \\
	&&	+ \frac{\tM}{2\alpha c}\sumi\int_{\T} \left|\partial_x \left(\ffi\right)\right|^2 \psixiinf \ dx 
		+ 4\alpha \tR^2E_m.  
\end{eqnarray*}
\noindent Finally, by grouping terms together and using the fact that $\alpha$ may be chosen such that $\tM/2\alpha c < 1$ (an appropriate choice for $\alpha$ is given later in the proof), we have
\begin{eqnarray}
\frac{1}{2}\frac{dP}{dt}
	&\leq& 	- \left(1 -\frac{\tM}{2\alpha c}\right)\sumi\int_{\T} \left|\partial_x \left(\ffi\right)\right|^2\psixiinf \ dx
		- \int_{\T} \lambda(x)\left|\psixi_0-\psixi_1\right|^2\frac{1}{\psixiinf} \ dx 
		+ 4\alpha \tR^2E_m \nonumber \\
	&\leq& 	- \left(1 -\frac{\tM}{2\alpha c}\right)\left[\sumi\int_{\T} \left|\partial_x \left(\ffi\right)\right|^2\psixiinf \ dx
		+ \int_{\T} \lambda(x)\left|\psixi_0-\psixi_1\right|^2\frac{1}{\psixiinf} \ dx \right]
		+ 4\alpha \tR^2E_m \nonumber \\
	&\leq& 	-\left(1 -\frac{\tM}{2\alpha c}\right)\sumi\int_{\T} (\psixi \L_i \psixi )
		\frac{1}{\psixi_{\infty}}\ dx
		+ 4\alpha \tR^2E_m \nonumber \\
	&\leq& 	-\left(1 -\frac{\tM}{2\alpha c}\right)\theta P
		+ 4\alpha \tR^2E_m,  \label{eq:firstP}
\end{eqnarray}
where the last line is a result of~\eqref{eq:spec_gap}, with $f:=\psixi - \psixi_{\infty}$.  Notice that $f \in \W_0$ since the normalization for $\psixi$ is $\displaystyle \sumi \int_{\T}\psixi dx = 1$.  

To complete the proof of Theorem \ref{th:main}, we now need to study the system of inequalities~\eqref{eq:dEm_dt} and~\eqref{eq:firstP}.
\subsection{Completing the proof}
To show that $E_m$ decays exponentially fast, we study the system of two inequalities \eqref{eq:dEm_dt} and \eqref{eq:firstP}.  Since, from \eqref{eq:E_c}, $E_c \leq P$, the system to be studied is  
\begin{equation}
\label{eq:sys_ineq}
\left\{
\begin{array}{l l}
\displaystyle\frac{dE_m}{dt} &\leq -2\left(\rho -\tR^2\varepsilon \right) E_m  
			+2\rho P
			+\displaystyle \frac{1}{4\varepsilon}F_0 \,\text{e}^{-8\pi^2t},  \\[1.5ex]
\displaystyle\frac{dP}{dt}&\leq 8\alpha \tR^2E_m
			-2\left(1 -\displaystyle\frac{\tM}{2\alpha c}\right)\theta P.  
\end{array}
\right.
\end{equation}
The parameters $\alpha>\frac{\tM}{2c}$ and $\varepsilon > 0$ remain to be chosen in order to obtain an exponential convergence with the best possible rate.  To fix $\alpha$, let us first study the eigenvalues of the matrix of coefficients, neglecting the terms in $\varepsilon$.
\begin{lemma}
Let us assume [H4]. The matrix
\begin{equation*}
A= \left(
\begin{array}{c c}
 -\rho & \rho \\ 4\alpha R^2 & -\left(1-\frac{\tM}{2\alpha c}\right)\theta
\end{array}
\right)  
\end{equation*}
is negative definite and $\alpha$ may be chosen so that the eigenvalues $-\lambda_{\pm}$ of $A$ are such that 
$$-\lambda_- \le -\lambda_+ = - \Lambda(\theta) < 0$$ 
where $\Lambda: (\theta_{\rm min},\infty) \rightarrow (0,\rho)$ is a positive, increasing function. We recall that $\displaystyle \theta_{\rm min}=\frac{8 \tM \tR^2}{c}$, where $\tR= C+M\rho^{-1/2}$. The function $\Lambda$ is such that $\Lambda \rightarrow 0$ as $\theta \rightarrow \theta_{\rm min}$ and $\Lambda \rightarrow \rho$ as $\theta \rightarrow \infty$. Moreover, $\displaystyle\Lambda( \rho + 2 \theta_{\rm min}) = \frac{\rho}{2}$.
\end{lemma}
\proof In order to prove the negative definiteness of the matrix $A$, we show that for certain values of $\alpha>0$, $\text{tr}(A)<0$ and $\text{det}(A)>0$. In the following, we only consider positive values of $\alpha$ (which is imposed by the previous computations). We have 
\begin{equation}
\label{eq:alpha}
 \text{tr}(A) = -\rho - \left(1-\frac{\tM}{2\alpha c}\right)\theta < 0
    ~~\text{iff}~~
 \alpha > \frac{\tM \theta}{2c(\rho+\theta)}  
\end{equation}
and 
\begin{equation}
\label{eq:detalpha}
 \text{det}(A) = \theta\rho\left(1-\frac{\tM}{2\alpha c}\right)-4\alpha R^2 \rho > 0
    ~~\text{iff}~~
 \alpha \in (\alpha_-,\alpha_+), ~
	\alpha_{\pm} = \frac{\theta c \pm \sqrt{\theta^2c^2 - 8\tM \theta R^2c}}{8R^2 c}.  
\end{equation}
The interval $(\alpha_-,\alpha_+)$ is indeed well defined and included in $[0, \infty)$ since $\theta> \theta_{\rm min}= 8\tM R^2 / c$ (hypothesis [H4]).  We seek an optimal $\alpha$ that minimizes eigenvalue $-\lambda_+$ and satisfies~\eqref{eq:alpha} and~\eqref{eq:detalpha}.  An analytical solution cannot be easily obtained.  We choose  
$$ \alpha = \alpha^* := \frac{\tM}{c}, $$
which appears to be very close to the optimal choice, from numerical computations. Notice that $\alpha^*$
 satisfies~\eqref{eq:alpha} and~\eqref{eq:detalpha}  since, for $\alpha=\alpha^*$, $\text{tr}(A)= -\rho - \theta/2< 0$ and $\text{det}(A) = \theta \rho / 2 - 4 R^2 \tM \rho / c >0$.  The eigenvalues of the matrix are now given by
\begin{equation*}
 -\lambda_{\pm}= \frac{1}{2} \left( -\left(\rho + \frac{\theta}{2}\right) 
		\pm\sqrt{\left(\rho - \frac{\theta}{2}\right)^2 + \frac{16 R^2 \tM \rho}{c}} \right) ~~< 0.  
\end{equation*}
The rate of convergence of the system is given by the largest of the two eigenvalues $-\lambda_+$. Let us introduce the function 
$$\Lambda(\theta)=-\frac{1}{2} \left( -\left(\rho + \frac{\theta}{2}\right) 
		+\sqrt{\left(\rho - \frac{\theta}{2}\right)^2 + \frac{16 R^2 \tM \rho}{c}} \right)$$
such that $\lambda_+=\Lambda(\theta)$.  It is easily shown that $\Lambda$ is an increasing function of $\theta$ with 
\begin{equation*}
\Lambda(\theta) \rightarrow 
\left\{
\begin{array}{r l}
 0      &\text{as}~~\theta \rightarrow \frac{8\tM R^2}{c},\\
 \rho &\text{as}~~\theta \rightarrow \infty.
\end{array}
\right.
\end{equation*}
Moreover, it is easy to check that $\displaystyle\Lambda( \rho + 2 \theta_{\rm min}) = \frac{\rho}{2}$, which concludes the proof.
 \endproof

We are now in position to complete the proof of Theorem~\ref{th:main}. Let us define $\Y(t) = (E_m(t),P(t))$.  Using~\eqref{eq:sys_ineq} the fact that $E_m~\leq~\|\Y\|_2$ (where $\|Y\|_2$ denotes the Euclidean norm of the two-dimensional vector $Y$), we obtain
\begin{eqnarray*}
  \frac{1}{2}\frac{d}{dt}\|\Y\|_2^2
	&=& \frac{1}{2}\frac{d}{dt}(E_m^2+P^2)\\
	&\le &     2\Y^TA\Y + 2R^2\varepsilon E_m^2 + \frac{1}{4\varepsilon}F_0 \,\text{e}^{-8\pi^2t}E_m\\
	&\leq&  -2 \Lambda(\theta) \, \|\Y\|_2^2 + 2R^2\varepsilon \|\Y\|_2^2 + \frac{1}{4\varepsilon}F_0\,\text{e}^{-8\pi^2t}\|\Y\|_2,
\end{eqnarray*}
\noindent and as a result,
\begin{eqnarray}
\label{eq:dV}
  \frac{d\|\Y\|_2}{dt}
	&\leq&   -2(\Lambda(\theta) - R^2\varepsilon)\|\Y\|_2 + \frac{1}{4\varepsilon}F_0\,\text{e}^{-8\pi^2t}.  
\end{eqnarray}
For arbitrary small $\varepsilon > 0$, let us consider $\lameps = \Lambda(\theta) - R^2\varepsilon < \Lambda(\theta)$. We may assume without loss of generality that $\lameps\neq 4\pi^2$. Then, from~\eqref{eq:dV}, one gets:
\begin{equation}
\label{eq:Em_conv_result}
 E_m\leq \|\Y\|_2
	\leq  K_\varepsilon \,
			\text{e}^{-2\min\{\lameps,4\pi^2\}t},
\end{equation}
where 
$$ K_\varepsilon = 2~\max\left\{\sqrt{E_m^2(0)+P^2(0)},
				\frac{F_0}{8\varepsilon|\lameps-4\pi^2|}\right\},$$
which concludes the proof of~\eqref{eq:CVEm}.

The exponential convergence of the total entropy $E$ results from the relation $E=E_M+E_m$, ~\eqref{eq:Em_conv_result} and Lemmas~\ref{lem:macro} and~\ref{lem:macro_fisher}.  The Csiszar-Kullback inequality implies the same for $\|\psi(t,\cdot)-\psi_{\infty}\|^2_{L^1}$.  

Finally, the convergence results on $A_t'$ are easily obtained from Lemma~\ref{lem:At_A_Em} and the fact that $\psixibar$ is bounded from below by a positive constant for times larger than an arbitrary small positive time, see the beginning of Section~3.3.2 in~\cite{lrs_nonlinearity} for more details.

\section*{Acknowledgements}
   This work is supported in part by the Agence Nationale de la Recherche, under the grant MEGAS (ANR-09-BLAN-0216-01).  The authors are grateful to Benjamin Jourdain, Eric Canc\`es, Andreas Eberle, Felix Otto and Gabriel Stoltz for helpful discussions.

\bibliographystyle{plain}

\begin{thebibliography}{10}
\bibitem{logsob2000}
C.~An\'e, S.~Blach\`ere, D.~Chafa\"{i}, P.~Foug\`eres, I.~Gentil, F.~Malrieu,
  C.~Roberto, and G.~Scheffer.
\newblock {\em Sur les in\'egalit\'es de {S}obolev logarithmiques}.
\newblock SMF, 2000.

\bibitem{arnold2001}
A.~Arnold, P.~Markowich, G.~Toscani, and A.~Unterreiter.
\newblock On logarithmic {S}obolev inequalities and the rate of convergence to
  equilibrium for {F}okker-{P}lanck type equations.
\newblock {\em Comm. Partial Differential Equations}, 26:35--43, 2001.

\bibitem{bakry1984}
D.~Bakry and M.~Emery.
\newblock Hypercontractivit\'e de semi-groupes de diffusion.
\newblock {\em C. R. Acad. Sci Paris S\'er. I}, 299:775--778, 1984.

\bibitem{ciccotti2008}
G.~Ciccotti, T.~Leli\`evre, and E.~Vanden-Eijnden.
\newblock Projection of diffusions on submanifolds: Application to mean force
  computation.
\newblock {\em Comm. Pure Appl. Math.}, 61:3, 2008.

\bibitem{darve-pohorille-01}
E.~Darve and A.~Pohorille.
\newblock Calculating free energies using average forces.
\newblock {\em J. Chem. Phys.}, 115(20):9196--9183, 2001.

\bibitem{otter1998}
W.K. den Otter and W.J. Briels.
\newblock The calculation of free energy differences by constrained molecular
  dynamics simulations.
\newblock {\em J. Chem. Phys.}, 109:4139, 1998.

\bibitem{gross1975}
L.~Gross.
\newblock Logarithmic {S}obolev inequalities.
\newblock {\em American Journal of Mathematics}, 97:1061--1083, 1975.

\bibitem{henin2004}
J.~H\'enin and C.~Chipot.
\newblock Overcoming free energy barriers using unconstrained molecular
  dynamics simulations.
\newblock {\em J. Chem. Phys.}, 121:2904, 2004.

\bibitem{henin2010}
J.~H\'enin, G.~Fiorin, C.~Chipot, and M.~L. Klein.
\newblock Exploring multidimensional free energy landscapes using
  time-dependent biases on collective variables.
\newblock {\em J. Chem. Theory Comput.}, 6:35--47, 2010.

\bibitem{holleystroock}
R.~Holley and D.~Stroock.
\newblock Logarithmic {S}obolev inequalities and stochastic {I}sing models.
\newblock {\em J. Stat. Phys.}, 46:1159--1194, 1987.

\bibitem{jlr_2009}
B.~Jourdain, T.~Leli\`evre, and R.~Roux.
\newblock Existence, uniqueness and convergence of a particle approximation for
  the adaptive biasing force process.
\newblock {\em to appear in ESAIM: M2AN}, 2010.

\bibitem{lelievre2010}
T.~Leli\`evre.
\newblock A general two-scale criteria for logarithmic {S}obolev inequalities.
\newblock {\em J.Func. Anal.}, 256(7):2211--2221, 2009.

\bibitem{lrs2007}
T.~Leli\`evre, M.~Rousset, and G.~Stoltz.
\newblock Computation of free energy profiles with parallel adaptive dynamics.
\newblock {\em J. Chem. Phys.}, 126:134111, 2007.

\bibitem{lrs_nonlinearity}
T.~Leli\`evre, M.~Rousset, and G.~Stoltz.
\newblock Long-time convergence of an adaptive biasing force method.
\newblock {\em Nonlinearity}, 21(6):1155--1181, 2008.

\bibitem{minoukadeh2010}
K.~Minoukadeh, C.~Chipot, and T.~Leli\`evre.
\newblock Potential of mean force calculations: A multiple-walker {A}daptive
  {B}iasing {F}orce approach.
\newblock {\em J. Chem. Theory Comput.}, 6:1008--1017, 2010.

\bibitem{otto2000}
F.~Otto and C.~Villani.
\newblock Generalization of an inequality by {T}alagrand, viewed as a
  consequence of the logarithmic {S}obolev inequality.
\newblock {\em J. Funct. Anal.}, 173(2):361--400, 2000.

\bibitem{sprik1998}
M.~Sprik and G.~Cicotti.
\newblock Free energy from constrained molecular dynamics.
\newblock {\em J. Chem. Phys.}, 109:7737--7744, 1998.

\bibitem{villanitopics}
C.~Villani.
\newblock {\em Topics in {O}ptimal {T}ransportation}.
\newblock American Mathematical Society, 2003.

\end{thebibliography}

\end{document}